\documentclass[12pt,fleqn]{amsart}
\usepackage[mathscr]{eucal}
\usepackage{amssymb,amsmath}
\usepackage{amsfonts}
\usepackage{a4}
\usepackage{xcolor}
\usepackage{enumerate}
\usepackage{mathtools}
\usepackage{lipsum}
\usepackage{tikz-cd}
\usepackage{geometry}
\usepackage{layout}
\usepackage[geometry]{ifsym}
\theoremstyle{plain}
\newtheorem{thm}{Theorem}[section]
\newtheorem{dfn}[thm]{Definition}
\newtheorem{prp}[thm]{Proposition}
\newtheorem{cor}[thm]{Corollary}
\newtheorem{lem}[thm]{Lemma}

\newtheorem{ex}[thm]{Example}

\theoremstyle{definition}
\newtheorem{definition}{Definition}[section]
\theoremstyle{remark}
\newtheorem{prob}[thm]{Problem}
\newtheorem{rem}[thm]{Remark}

\numberwithin{equation}{section}

\newcommand{\FF}{\mathcal{F}}

\newcommand{\eps}{\varepsilon}

\DeclareMathOperator{\supp}{supp}

\DeclareMathOperator{\wnt}{int}

\newcommand{\phe}{\varphi}

\newcommand{\eco}{\mathcal{EC}_{\omega c}}
\newcommand{\ny}{\mathcal{NY}}
\newcommand{\elk}{\mathcal{L}_\kappa}

\newcommand{\ult}{{\rm ult}}

\newsavebox{\overlongequation}

\newcommand{\cont}{\mathfrak c}

\newcommand{\fB}{\mathfrak B}

\newcommand{\cA}{\mathcal A}
\newcommand{\cI}{\mathcal I}
\newcommand{\cC}{{\mathcal C}}
\newcommand{\cF}{{\mathcal F}}

\newcommand{\cU}{{\mathcal U}}
\newcommand{\cV}{{\mathcal V}}

\newcommand{\er}{\mathbb R}
\newcommand{\qu}{\mathbb Q}
\newcommand{\vf}{\varphi}

\newcommand{\stevo}{Todor\v{c}evi\'c}

\newcommand{\sm}{\setminus}

\newcommand{\en}{\mathbb N}
\newcommand{\sub}{\subseteq}

\newcommand{\ord}{\protect{\rm ord} }

\newcommand{\sikg}{\protect{ \Sigma_{\kappa}([0,1]^{\Gamma}) }}
\newcommand{\sikl}{\protect{ \Sigma_{\kappa}([0,1]^{\lambda}) }}
\newcommand{\Lra}{\Leftrightarrow}
\DeclareMathOperator\BTU{\mbox{\BigTriangleUp}}

\begin{document}

\title{Digging into the classes of $\kappa$-Corson compact spaces}

\author[W.\ Marciszewski]{Witold Marciszewski}
\address{Institute of Mathematics\\
University of Warsaw\\ Banacha 2\newline 02--097 Warszawa\\
Poland\\
ORCID identifier: 0000-0003-3384-5782}
\email{wmarcisz@mimuw.edu.pl}

\author[G.\ Plebanek]{Grzegorz Plebanek}
\address{Instytut Matematyczny\\ Uniwersytet Wroc\l awski\\ Pl.\ Grunwaldzki 2\\
50-384 Wroc\-\l aw\\ Poland} \email{grzegorz.plebanek@math.uni.wroc.pl}

\author[K.\ Zakrzewski]{Krzysztof Zakrzewski}
\address{Institute of Mathematics\\
University of Warsaw\\ Banacha 2\newline 02--097 Warszawa\\
Poland\\}
\email{Krzysztof.Zakrzewski@mimuw.edu.pl}

\thanks{The first and the third authors were partially supported by the NCN
(National Science Centre, Poland) research grant no.\ 2020/37/B/ST1/02613.
The second author was supported by the NCN
(National Science Centre, Poland), under the Weave-UNISONO call in the Weave
programme 2021/03/Y/ST1/00124.
}

\date{\today}

\begin{abstract}
For any cardinal number $\kappa$ and an index set $\Gamma$,
the $\Sigma_\kappa$--product of real lines  $\Sigma_\kappa(\mathbb{R}^\Gamma)$
consists of all elements of $\er^\Gamma$ having $<\kappa$ many nonzero coordinates.
A compact space $K$ is $\kappa$-Corson compact if it can be embedded into
$\Sigma_\kappa(\mathbb{R}^\Gamma)$ for some $\Gamma$.

The class of ($\omega_1$)-Corson compact spaces has been intensively studied over last decades.
We discuss  properties of $\kappa$-Corson compacta for various cardinal numbers $\kappa$
as well as properties of related Boolean algebras and spaces of continuous functions.

We present here a detailed discussion of the class of $\omega$-Corson compacta
extending the results of Nakhmanson and Yakovlev \cite{NY}.
For $\kappa>\omega$, our results on $\kappa$-Corson compact spaces are related to the
line of research  originated by Kalenda \cite{Ka} and Bell and Marciszewski \cite{BM04},
and continued by Bonnet, Kubi\'s and \stevo\  in their  recent paper \cite{BKT}.
\end{abstract}

\maketitle


\section{Introduction}\label{intro}

Given an infinite cardinal number $\kappa$, a compact space $K$ is
said to be {\em $\kappa$-Corson compact} if, for some set $\Gamma$,
$K$ is homeomorphic to a subset of the $\Sigma_\kappa$--product of real lines, that is
\[\Sigma_\kappa(\mathbb{R}^\Gamma)=\{x\in \mathbb{R}^\Gamma: |\{\gamma: x(\gamma)\neq 0\}| < \kappa\}.\]
For $\kappa=\omega_1$ this definition names the well-studied class of Corson compacta,
compact spaces that can be embedded into the familiar $\Sigma$-product of real lines, see e.g.\ \cite{Ne} and \cite{AMN88}.
For arbitrary $\kappa$,  $\kappa$-Corson compact spaces were first considered by Kalenda \cite{Ka}
and in \cite{BM04} (where the terminology was slightly different).

The classes of $\kappa$-Corson compacta are a subject  of a recent paper by Bonnet, Kubi\'s and
\stevo\ \cite{BKT} which motivated us to undertake  our research that we report here.
Trying to examine them systematically one discovers quickly that the case $\kappa=\omega$ is
quite special. In sections \ref{omega_corson} -- \ref{sec_eco} we investigate  $\omega$-Corson compacta
and a wider, somewhat more natural,  class of compact spaces that can be embedded into
a $\sigma$-product of metrizable compacta. The latter
we call  NY compacta to stress
the contribution of  Nakhmanson and Yakovlev \cite{NY} to the subject.

For regular $\kappa>\omega$ there is a number of properties of $\kappa$-Corson compacta that
may be seen as natural analogues of known features of the usual Corson compact spaces.
However, if $\kappa$ is a singular cardinal, then one faces number of problems that we were not able to
handle, with the exception of results from section \ref{local_properties}.

Here is the summary of our main results:

\begin{enumerate}[(A)]
\item Extending results from \cite{NY}, we give internal characterizations
of $\omega$-Corson compacta (Theorem \ref{charact_omega_Corson})  and of  NY compacta (Theorem \ref{charact_ec_sigma_m}). We discuss the stability of these classes of compacta under countable products (Corollary \ref{cor_prod_not_ny}), continuous images (Corollary \ref{cor:ny_cont_im}), and the functor $P$ which assigns to a compact space $K$ the compact space $P(K)$
of regular probability Borel measures equipped with the $weak^\ast$ topology (Remark \ref{conv_no_NY}). In \cite{NY} it was proved that all NY compacta are hereditarily metacompact, in section \ref{sec_more_prop} we investigate this covering property for Eberlein compacta not belonging to the class $\ny$ (Example \ref{ex_Grunehage}, Theorem \ref{thm_prod_not_HM}).
\item We present some results on $\pi$-character and tightness of  $\kappa$-Corson
compact space for arbitrary $\kappa>\omega$, answering two question of Kalenda \cite{Ka}.
\item We prove  that for every uncountable regular $\kappa$ the class
of $\kappa$-Corson compact spaces is stable under continuous images (Theorem \ref{stable}).
The result was first proved in \cite{BM04} for succesor cardinals and in \cite{BKT} for
arbitrary regular $\kappa$ by some model-theoretic approach. Our proof extends
the argument of Gul'ko \cite{Gu77} for $\kappa=\omega_1$.
\item In section \ref{bds} we partially solve a problem posed in \cite{BKT}:
Theorem \ref{main} says, in particular, that no $\sigma$-complete Boolean algebra of regular cardinality
$\kappa$ can have a $\kappa$-Corson compact Stone space.
\item In section \ref{sec_CpK} we discuss   generalizations of a
result due to Pol  \cite{Po84},
characterizing  Corson compacta
$K$ in terms of topological properties of function spaces  $C_p(K)$. We extend this characterization for $\kappa$-Corson compacta, where $\kappa$ is regular (Theorem \ref{charact_k_Corson_CpK}), and use it to present yet another  proof of the stability of the classes of these compacta under continuous images.
\item In section \ref{measures}  we show that the question for
which $\kappa$ the class of $\kappa$-Corson compact spaces is closed under the mentioned above functor $P$ (of
forming the space of probability measures) is naturally related to the notion of  a caliber of measures
(so is, typically, undecidable in the usual set theory). In turn, in the final section we discuss
the stability of $\kappa$-Corson compacta with respect to Banach space isomorphisms.
\end{enumerate}

The authors are very grateful to the referee for a very careful reading and several comments that enabled us to improve the presentation.

\section{Preliminaries}\label{prelim}

\subsection{Terminology and notation}
All topological spaces under consideration are assumed to be completely regular (Tikhonov).
As usual, the topological weight of a space $X$ is denoted by $w(X)$
(see section \ref{local_properties} for other cardinal invariants).

For a set $X$ and $n\in\omega$, we use the standard notation
\[ [X]^n=\{A\subseteq X: |A|=n\}, \quad [X]^{\le n}=\bigcup\{[X]^k: k\le n\};\]
moreover $[X]^{< \omega}$ stands for $\bigcup\{[X]^k: k< \omega\}$.

A family $\mathcal{U}$ of subsets of a space $X$
is $T_0$-separating if, for every pair of distinct points $x,y$ of
$X$, there is $U\in\mathcal{U}$ containing exactly one of the
points $x,y$.

Given a family $\mathcal{U}$ of subsets of a space $X$ we define $\ord(x,\mathcal{U}) = |\{U\in \mathcal{U}: x\in U\}|$ for $x\in X$. We say that $\mathcal{U}$ is \emph{point-finite} if $\ord(x,\mathcal{U}) < \omega$ for all $x\in X$.

For a locally compact space $X$, $\alpha(X)$ denotes the one point compactification of $X$. We denote the point at infinity of this compactification by $\infty_X$. For an infinite cardinal number $\kappa$, $D(\kappa)$ denotes a discrete space of cardinality $\kappa$, and $A(\kappa) = \alpha(D(\kappa))$.

Recall that for any compact space $K$,  its  Aleksandrov duplicate $AD(K)$ is defined as
the space $AD(K) = K\times 2$ in which all the points $(x,1)$  for $x\in K$ are isolated
and basic neighborhoods of a point $(x,0)$ have the form $(U\times 2)\setminus \{(x,1)\}$, where $U$ is an open neighborhood of $x$ in $K$.

A space $X$ is \textit{scattered} if no nonempty subset $A\subseteq X$ is dense-in-itself.
For a scattered space $K$, its
\textit{Cantor-Ben\-dixson height} $ht(X)$  is the minimal ordinal number
$\alpha$ such that the Cantor-Bendixson derivative $K^{(\alpha)}$
of the space $K$ is empty. The Cantor-Ben\-dixson height of a
compact scattered space is always a nonlimit ordinal.

A space $X$ is \textit{metacompact} if every open cover $\mathcal{U}$ of $X$ has a point-finite open refinement $\mathcal{V}$. It is \textit{$\sigma$-metacompact} if every open cover $\mathcal{U}$ of $X$ has an open refinement $\mathcal{V} = \bigcup_{n\in\omega}\mathcal{V}_n$, where each family $\mathcal{V}_n$ is point-finite.

For a space $X$, the sets of the form $f^{-1}((0,1])$, where $f:X \to [0,1]$ is a continuous function, are called the \textit{cozero sets}.

Recall that a space $X$ is  \emph{strongly countable-dimensional} if $X$ is a countable union of closed subspaces of finite covering dimension, see \cite{En}.

\subsection{$\Sigma_\kappa$-products}\label{subsec_sigma_prod}
We often write $I$ for  the closed unit interval $[0,1]$.
In the sequel, $\kappa$ always stands for an infinite cardinal number.

Given $x$ in the product space $\er^\Gamma$, we write
$\supp(x)=\{\gamma\in\Gamma: x(\gamma)\neq 0\}$.
Thus $\Sigma_\kappa(\er^\Gamma)$ is the space of all $x\in\er^\Gamma$ satisfying $|\supp(x)|<\kappa$
and $\Sigma_\kappa ([0,1]^\Gamma)$ is defined analogously.
Traditionally, $\Sigma_{\omega_1}(\er^\Gamma)$ is written as $\Sigma(\er^\Gamma)$
and $\Sigma_{\omega}(\er^\Gamma)$ is written as $\sigma(\er^\Gamma)$.
However, we need to specify more notation for the case $\kappa=\omega$.

Let $\{X_\gamma: \gamma\in \Gamma\}$ be a family of nonempty topological spaces, and let $a_\gamma$ be
a fixed point in $X_\gamma$. The $\sigma$-product of the family
$\{(X_\gamma,a_\gamma): \gamma\in \Gamma\}$ is defined as
\[\sigma(X_\gamma,a_\gamma,\Gamma) = \{(x_\gamma)_{\gamma\in \Gamma}\in \prod_{\gamma\in \Gamma} X_\gamma: |\{\gamma\in\Gamma: x_\gamma \ne a_\gamma\}| < \omega\}\,.\]

If $X_\gamma = I$ and $a_\gamma = 0$, for all $\gamma\in \Gamma$,  then
$\sigma(X_\gamma,a_\gamma,\Gamma)$ is simply denoted by $\sigma(I,\Gamma)$.
Likewise,
if $X_\gamma = I^\omega$ and $a_\gamma = (0,0,\dots)$, for all $\gamma\in \Gamma$,
then  $\sigma(X_\gamma,a_\gamma,\Gamma)$ is denoted by $\sigma(I^\omega,\Gamma)$.

\subsection{Corson ($\kappa$-Corson) compacta and subclasses}\label{subsec_kappa_corson}

As we have already mentioned, a compact space is said to be $\kappa$-Corson compact if it can be embedded into
$\Sigma_\kappa(\er^\Gamma)$ for some $\Gamma$.
To respect the existing tradition concerning the case $\kappa=\omega_1$, we write Corson compact rather than
$\omega_1$-Corson compact.

\begin{rem}\label{basic_observations}
We record here the following preliminary  observations.
\begin{enumerate}[(a)]
\item While defining $\kappa$-Corson compacta,   one can clearly replace  $\Sigma_\kappa(\er^\Gamma)$ by
$\Sigma_\kappa([0,1]^\Gamma)$.
\item If $K$ is $\kappa$-Corson compact, then $K$ embeds into
$\Sigma_\kappa([0,1]^\Gamma)$ where $\Gamma$ is of size $w(K)$.
\item {\em Every} compact space $K$ is $\kappa$-Corson compact for $\kappa> w(K)$.
\item A compact space $K$ is $\kappa$-Corson compact if and only if
there is a family $\cF$ of continuous functions $K\to [0,1]$ that separates the points of $K$ and
$|\{f\in\cF: f(x)\neq 0\}|<\kappa$ for every $x\in K$.
\end{enumerate}
\end{rem}

The following result is a straightforward generalization of the well known
Rosenthal-type characterization of Corson compacta (see e.g.\  \cite{BKT}).

\begin{prp}\label{charact_BKT}
Let $\kappa$ be an uncountable cardinal number. For a compact space $K$, the following
conditions are equivalent
\begin{enumerate}[(a)]
\item $K$ is $\kappa$-Corson;
\item There exists a family $\mathcal{U}$ consisting of cozero subsets of $K$ which is $T_0$-separating, and
$\ord(x,\mathcal{U}) < \kappa$ for all $x\in K$.	
\end{enumerate}
\end{prp}

A space $K$ is an \emph{Eberlein} compact
if $K$ is homeomorphic to a weakly compact subset of a Banach
space.
By the celebrated Amir-Lindenstrauss theorem, a compact space $K$ is an Eberlein compactum if and only if
it can be embedded into
\[c_0(\Gamma)=\{x\in \mathbb{R}^\Gamma:
 \big| \{\gamma: |x(\gamma)|>\varepsilon\}\big| <\omega
\mbox{ for every } \varepsilon>0 \},\]
for some set $\Gamma$, see \cite{Ne}.
As $c_0(\Gamma)\sub \Sigma_{\omega_1}(\er^\Gamma)$,
the class of Corson compact spaces contains all Eberlein compacta.

If $K$ is homeomorphic to a weakly compact subset of a
Hilbert space, then $K$ is said to be  \emph{uniform Eberlein
compact}. Note that the class of uniform Eberlein compacta contains all
metrizable compact spaces.

\subsection {Boolean algebras}\label{subsection_ba}
For a Boolean algebra $\fB$ we denote by $\ult(\fB)$ its Stone space (of all ultrafilters on $\fB$);
the Stone isomorphism
\[\fB\ni b\mapsto\widehat{b}=\{p\in\ult(\fB): b\in p\}\]
identifies $\fB$ with the algebra of clopen subsets of $\ult(\fB)$.

We say that  a Boolean  algebra $\fB$ is {\em $\kappa$-Corson}  if
its Stone space $\ult(\fB)$ is $\kappa$-Corson compact. The property can be characterized internally
as follows.

\begin{lem}\label{lemma:ba}
A Boolean  algebra $\fB$ is $\kappa$-Corson if and only if there is
$G\sub \fB$ that generates $\fB$ and every centered subfamily $G_0\sub G$ has size $<\kappa$.
\end{lem}

The proof of the lemma above can be found in \cite{BKT}; it is an obvious analogue
of the case $\kappa=\omega_1$ which was already noted in \cite[Lemma 2.2.]{MP17}.

\subsection{Function spaces}

Given a compact space $K$, by $C(K)$ we denote the Banach space of
continuous real-valued functions on $K$, equipped with the standard
supremum norm.
As usual, we identify the dual space $C(K)^\ast$ with the space
$M(K)$ of signed Radon measures on $K$ having finite variation
and write ${\mu}(g)$ for $\int_K g\;{\rm d}\mu$. The symbol $M_1(K)$ stands for the unit ball of
$M(K)$, equipped with the $weak^\ast$ topology inherited from $C(K)^\ast$.
Finally, $P(K)$ is the subset of regular probability measures on $K$.

We will also consider   function spaces of the form $C_p(K)$; in such a case   $C(K)$ is equipped with the pointwise
convergence topology.

\section{$\omega$-Corson compact spaces and NY compacta}\label{omega_corson}

Motivated by the results by Nakhmanson and  Yakovlev  \cite{NY}, we introduce here
the following notion.

\begin{dfn}\label{ny_def}
We shall say that a compact space $K$ is {\em NY compact}
or write $K\in \ny$
 if $K$
can be embedded into some $\sigma$-product of metrizable  compacta.
\end{dfn}

We start by the following easy observations.

\begin{prp}\label{prp_sigma_prod} Let $K$ be a  compact space $K$.
\begin{enumerate}[(a)]
\item $K$ is $\omega$-Corson if and only if it can be embedded into some
$\sigma$-product of metrizable finitely dimensional compacta.
\item $K$ is NY compact if and only if it can be embedded into the $\sigma$-product $\sigma(I^\omega,\Gamma)$ for some set $\Gamma$.
\end{enumerate}	
\end{prp}

\begin{proof}
The forward implication in (a) follows from Remark \ref{basic_observations}(a).
To check the reverse implication, suppose that
  $\Xi$ is an embedding of  $K$ into $\sigma(X_{\gamma},a_{\gamma},\Gamma)$, where $X_{\gamma}$'s
are compact metrizable spaces of finite dimension.

For every $\gamma\in\Gamma$ there exist a natural number $k_{\gamma}$ and
a homeomorphic embedding $\Phi_{\gamma}:X_{\gamma}\hookrightarrow [-1,1]^{k_{\gamma}}$
such that $\Phi_{\gamma}(a_{\gamma})=\mathbf{0}^{k_{\gamma}}$; hence we can treat $K$ as a subspace of $\sigma([-1,1],0,{\Gamma'})$ for a suitable set $\Gamma'$.
We identify $[-1,1]$ with $(\{0\} \times [0,1]) \cup ([0,1] \times \{0\})$ by rotating $[-1,1]$ over 90 degrees. In this way we turn a subspace of $\sigma([-1,1],0,{\Gamma'})$ into a subspace of $\Sigma_\omega([0,1]^{\Gamma'\times 2})$.

Part (b) follows immediately from the fact that that any metrizable compactum can be embedded into the Hilbert cube $I^\omega$ and $I^\omega$ is homogeneous.
\end{proof}

Let us note that Proposition \ref{prp_sigma_prod}(b) gives another equivalent definition of the 
class $\ny$ which will be often used below.

\begin{prp}\label{prp_omega_SCD}
Every $\omega$-Corson compact space is Eberlein compact and strongly countably dimensional.
\end{prp}

\begin{proof}
If $K$ is a $\omega$-Corson compact space, then it is also Eberlein compact by the very definition.
To verify that second statement we can assume that  $K\sub \sigma(I,\Gamma)$.

Note that
$\sigma(I,\Gamma)=\bigcup_n A_n$ where
\[A_{n}=\{x\in\sigma(I,\Gamma): |\supp(x)|\leq n\}.\]
Moreover, every $A_n$ is a compact space of dimension $n$, see \cite{EP}.
 It follows that  $\dim(A_n\cap K)\le \dim A_n = n$ and hence $K$ is strongly countably dimensional.
\end{proof}

\begin{prp}\label{prp_metr_omega_SCD}
A metrizable compact space $K$ is $\omega$-Corson if and only if it is strongly countably dimensional.
\end{prp}

\begin{proof}
This  follows immediately from Proposition \ref{prp_omega_SCD}
and the result from \cite{En} stating  that the space $\sigma(I,\omega)$
is universal  for the class of separable metrizable spaces that
are strongly countably dimensional.
\end{proof}

It becomes clear  that every $\omega$-Corson compact space is NY compact, while
the Hilbert cube is NY compact but not $\omega$-Corson compact.
We shall see later that  there are  zero-dimensional Eberlein compacta which are not
$\omega$-Corson.
The assertion of  Proposition \ref{prp_metr_omega_SCD} in fact
 holds  for NY compact spaces $K$, see Corollary \ref{prp_NY_omega_SCD}.

\begin{rem}\label{prp_boolean_omega_Corson}
Given a Boolean algebra $\fB$ and its Stone space $K=\ult(\fB)$, the following
conditions are equivalent:
\begin{enumerate}[(i)]
	\item $K$ is  $\omega$-Corson compact;
	\item $K$ has a $T_0$-separating, point-finite family $\mathcal{U}$ consisting of clopen subsets;
	\item $K$ is  scattered Eberlein compact;
	\item  $K$ is scattered Corson compact.
\end{enumerate}

Indeed, (i)$\Leftrightarrow$(ii) follows immediately from  Lemma \ref{lemma:ba} and
(ii) implies that $K$ can be embedded into some $\sigma(\{0,1\},\Gamma)$;
it is easy
to check that every compact subset of such a $\sigma$-product is scattered.
Then (iv)$\Rightarrow$(iii)  is a result of Alster \cite{Al79}.
\end{rem}

In connection with Remark \ref{prp_boolean_omega_Corson}, we shall indicate that
the Rosenthal-type characterization of
$\kappa$-Corson compacta given by Proposition \ref{charact_BKT}
fails for $\kappa=\omega$.

\begin{prp}\label{cos}
For a compact space $K$, the following
conditions are equivalent:
\begin{enumerate}[(a)]
		\item There exists a $T_0$-separating, point-finite family $\mathcal{U}$ consisting of cozero subsets of $K$;
		\item $K$ is a scattered Eberlein compact space.
\end{enumerate}
\end{prp}

In view of Remark  \ref{prp_boolean_omega_Corson}, Proposition \ref{cos}
follows immediately from the following lemma.

\begin{lem}
	 If a compact space $K$ has a point-finite, $T_{0}$ - separating family $\mathcal{U}$ of open sets,
then $K$ is scattered.
\end{lem}

\begin{proof}
Suppose  that $K$ is not scattered  and let   $A\subseteq K$ be a nonempty subset
without isolated points. To arrive at contradiction,
we inductively  construct  nonempty open sets $V_{n}$ such that $\overline{V_{n+1}}\subseteq V_{n}$ and $|\{U\in\mathcal{U}:V_{n}\subseteq U\}|\geq n$ for every $n$.

To start,  pick  $x,y\in A$, $x\neq y$.
There exists $V_{1}\in\mathcal{U}$ such that $|V_{1}\cap\{x,y\}|=1$;
say  $x\in V_{1}$. As $x$ is not isolated in $A$,
 there is $y_{1}\in A\cap V_{1}\sm\{x\}$. In turn,  there exists
$U_{1}\in\mathcal{U}$ separating $x$ and $y_{1}$.
Assume, for instance,  that $y_{1}\in U_{1}$ and  let $V_{2}$ be an open set such that
$y_{1}\in V_{2}\subseteq \overline{V_{2}}\subseteq U_{1}\cap V_{1}$.
Continuing in this fashion, we get the required sequences of $V_n$'s.

Now the set $F=\bigcap_{n} V_{n}=\bigcap_{n}\overline{V_{n}}$ is
nonempty by compactness. If $z\in F$ then $z$ belongs to infinitely many elements of
the family  $\mathcal{U}$, a contradiction.
 \end{proof}

We close this section by a few observations on the stability of our classes under
some operations. It was proved in \cite[Proposition 2.8]{Ma} that
the Aleksandrov duplicate $AD(K)$ of an (uniform) Eberlein compact space $K$ is (uniform) Eberlein compact.
Using the same argument  one  can show a counterpart of this observation for $\kappa$-Corson and NY compacta.

\begin{prp}\label{prp_AD_KC_NY} For any infinite cardinal number $\kappa$,
the Aleksandrov duplicate $AD(K)$ of a $\kappa$-Corson compact (NY compact) space $K$ is
$\kappa$-Corson compact (NY compact).
\end{prp}

\begin{prp}\label{prp_discrete_unions}
Let $\kappa>\omega $  and let $K_t,\ t\in T,$ be a family of
$\kappa$-Corson compact spaces (NY compact spaces).
Then the one point compactification $\alpha(\bigoplus_{t\in T} K_t)$ of the discrete union
$\bigoplus_{t\in T} K_t$  is $\kappa$-Corson compact  (NY compact, respectively).
\end{prp}

\begin{proof}
We will prove this fact only for $\kappa$-Corson compacta (for NY compacta the proof is very similar).

Denote the point at infinity of  $\alpha(\bigoplus_{t\in T} K_t)$ by $\infty$.
We can assume that $K_t\subseteq \Sigma_\kappa(\mathbb{R}^{\Gamma_t})$ for some set $\Gamma_t$, and the sets $\Gamma_t$ are pairwise disjoint and disjoint from $T$.

Consider
$\Gamma = T\cup \bigcup_{t\in T} \Gamma_t\,.$ Let $\psi: \alpha(\bigoplus_{t\in T} K_t) \to \Sigma_\kappa(\mathbb{R}^{\Gamma})$ be defined by
$$\psi(x)(\gamma) = \begin{cases} x(\gamma)& \text{if } x\in
K_t\,,\ \gamma\in \Gamma_t\,,\ t\in T \\
1& \text{if }  x\in
K_t\,,\ \gamma = t\,,\ t\in T \\
 0& \text{if } x = \infty\,, \gamma\in \Gamma
\end{cases}$$
for $x\in \alpha(\bigoplus_{t\in T} K_t)$. A routine verification shows that $\psi$ is an embedding.
\end{proof}

\section{Characterizations of $\omega$-Corson compacta and $NY$ compacta}\label{sec_charact}
\label{omega_corson2}

The main result of this section presents a characterization of $\omega$-Corson compact spaces
given in Theorem \ref{charact_omega_Corson}. Our theorem builds on a characterization of the class
$\ny$ obtained by Nakhmanson and Yakovlev  \cite{NY}, which we expand by
adding another condition, see Theorem \ref{charact_ec_sigma_m}(iv).

The following theorem offers several internal characterizations of compacta
from the class $\ny$.
Recall that a family $\mathcal{A}$ of subsets of a topological space $X$ is
closure preserving if $\overline{\bigcup \mathcal{A}'} = \bigcup\{\overline{A}: A\in \mathcal{A}'\}$
for any subfamily $\mathcal{A}'\subseteq \mathcal{A}$.


\begin{thm}\label{charact_ec_sigma_m}
For a compact space $K$, the following
conditions are equivalent:
\begin{enumerate}[(i)]
\item $K$ is NY compact;
\item there exists a $T_0$-separating family
$\mathcal{U} = \bigcup\{\mathcal{U}_\gamma: \gamma\in \Gamma\}$
consisting of cozero subsets of $K$, where each $\mathcal{U}_\gamma$ is countable and the family $\{\bigcup\mathcal{U}_\gamma: \gamma\in \Gamma\}$ is point-finite;
\item $K$ has a closure preserving cover consisting of metrizable compacta;
\item $K$ is hereditarily metacompact and each nonempty subspace $A$ of $K$
contains a nonempty relatively open subspace $U$ of countable weight.
\end{enumerate}
\end{thm}

The equivalences  (i)$\Lra$(ii)$\Lra$(iii) here were already proved in  \cite{NY}.
To incorporate condition (iv) we need some  preparations.
The  proposition given below was also proved in \cite[Theorem 11]{Ya}
using
the equivalent assumption that the compact space $K$ in question has a closure preserving cover
consisting of metrizable compacta. Our argument starting from the assumption that $K\in\ny$
seems to be shorter.

\begin{prp}\label{sigma-m otw zb ma otw z przel bazą}
If $K\in\ny$, then each nonempty subspace $A$ of $K$ contains a nonempty relatively open subspace $U$ of countable weight.
\end{prp}

\begin{proof}
By Proposition \refeq{prp_sigma_prod}(b), we can assume that $K\subseteq\sigma(I^{\omega},\Gamma)$ for some $\Gamma$. Obviously, it is enough to prove the assertion for the closure $\overline{A}$ of a nonempty set $A\subseteq K$. Since $\overline{A}$ is a compact subset of $\sigma(I^{\omega},\Gamma)$,
 for simplicity we can assume that $\overline{A} = K$.

Suppose  that there exists nonempty  open set $U\subseteq K$ such that
its every nonempty open subset has uncountable weight.
We can assume that  $U$ is a basic open set, that is
\[U=\{x\in K: x(\gamma)\in V_\gamma\mbox{ for } \gamma\in\Gamma_0\},\]
for  some finite $\Gamma_{0}\subseteq\Gamma$ and some family of open sets
$\{V_{\gamma}: \gamma\in \Gamma_0\}$ in $I^{\omega}$.

Note that  for every countable $\Gamma'\subseteq\Gamma$ there must be
$x\in U$ and $\gamma\in\Gamma\setminus\Gamma'$ such that
$x(\gamma)\neq 0$. Indeed, otherwise $U$ would be contained in a countable product of Hilbert's cubes
 so it would have a countable base.

Take $x\in U$ and $\gamma_{1}\in\Gamma\setminus\Gamma_{0}$ such that
$x(\gamma_{1})\neq 0$. Take open $G$ with  $x\in G\sub K$
and such that for every $y\in G$ we have   $y(\gamma_{1})\neq 0$.
Then find open  $W\subseteq K$  such that
$x\in W\subseteq \overline{W}\subseteq U$ and put  $U_{1}=W\cap G$.
Finally, set   $\Gamma_{1}=\Gamma_{0}\cup\{\gamma_{1}\}$.
Then for every $y\in U_{1}$ we have $y(\gamma_{1})\neq 0$ and
$x\in U_{1}\subseteq \overline{U_{1}}\subseteq U$.

Repeating the above argument, we can  inductively construct nonempty open subsets $U_{n}$ of $K$
and  distinct coordinates $\gamma_{n}\in\Gamma$ such that
$U_{n+1}\subseteq \overline{U_{n+1}}\subseteq U_{n}$ and
  if $z\in U_{n}$, then $z(\gamma_{n})\neq 0$. Clearly, by compactness,
$F=\bigcap_{n}U_{n}=\bigcap_{n}\overline{U_{n}}\neq\emptyset$ and
any $w\in F$ has infinitely many nonzero coordinates, a contradiction.
\end{proof}

For the next lemma we  slightly modify  the argument used in  \cite[Theorem 2]{NY}.

\begin{lem}
\label{pokr otw}
If $K\in\ny$ then for any family $\mathcal{V}$  of open subsets of $K$, there exists a point-finite  open refinement $\mathcal{U}$
such that $\bigcup \cU=\bigcup\cV$ and every $U\in {\mathcal U}$ is $\sigma$-compact.
\end{lem}

\begin{proof}

Consider    $K\subseteq \sigma(X_{\gamma},\xi_{\gamma},\Gamma)=Z$ where
$\{X_{\gamma}:\gamma\in \Gamma\}$ is a family of compact metrizable spaces and
$\xi_{\gamma}\in X_{\gamma}$ for every $\gamma\in \Gamma$.

Since, for a $\sigma$-compact $V\subseteq Z$, the intersection $V\cap K$ is also $\sigma$-compact it is enough to prove that for any family $\mathcal{V}$  of open subsets of $Z$, there exists a point-finite  open refinement $\mathcal{U}$
such that $\bigcup \cU=\bigcup\cV$ and every $U\in {\mathcal U}$ is $\sigma$-compact. Moreover, without loss of generality, we can assume that the family $\mathcal{V}$ consists of basic open sets in $Z$.

For $B\in [\Gamma]^{<\omega}$, let $\pi_{B}:Z\to \prod_{\gamma\in B}
X_{\gamma}$ be the projection restricted to $Z$. Observe that, for each open set $U \subseteq \prod_{\gamma\in B}
X_{\gamma}$, the set
\[\pi_{B}^{-1}(U) = U \times \sigma(X_{\gamma},\xi_{\gamma},\Gamma\setminus B)\]
is $\sigma$-compact being a product of
$\sigma$-compact spaces: a $\sigma$-product and
an open subset of a metrizable compact space. In particular, each basic open set in $Z$ is $\sigma$-compact, hence all $V\in \mathcal{V}$ are $\sigma$-compact.

Given $B\in [\Gamma]^{<\omega}$, let
\[U(B)=\{x\in Z:\forall_{\gamma\in B}\ x(\gamma)\neq \xi_{\gamma}\} = \pi_{B}^{-1}(\prod_{\gamma\in B}
(X_{\gamma}\setminus\{\xi_{\gamma}\})) .\]

Every $U(B)$ is open in $Z$ and $\sigma$-compact, and the family $\{U(B):B\in [\Gamma]^{<\omega}\}$ is point-finite.

For $B\in [\Gamma]^{<\omega}$, let
\[F(B)=\{x\in Z:x(\gamma)=\xi_{\gamma} \iff \gamma\in \Gamma\setminus B\}.\]
Observe that $\pi_{B}| F(B)$ maps $F(B)$ homeomorphically onto $\prod_{\gamma\in B}
(X_{\gamma}\setminus\{\xi_{\gamma}\})$.
By the definition of a $\sigma$-product, $\bigcup\{F(B):B\in [\Gamma]^{<\omega}\}=Z$.

Let  $\mathcal{K}(B)=\{V\cap F(B):V\in\mathcal{V}$\} 
for $B\in [\Gamma]^{<\omega}$. Since $F(B)$ is metrizable, there exists a point-finite family $\mathcal{M}(B)$ of open subsets of $F(B)$, such that $\mathcal{M}(B)\prec \mathcal{K}(B)$ and $\bigcup \mathcal{M}(B) = \bigcup  \mathcal{K}(B)$. Each set $W\in \mathcal{M}(B)$ is $\sigma$-compact because $F(B)$ is $\sigma$-compact and metrizable.

Let \[\mathcal{M}'(B)=\{\pi_{B}(W):W\in \mathcal{M}(B)\}.\]
Since $\pi_{B}| F(B)$ is a homeomorphism,
$\mathcal{M}'(B)$ is a point-finite family of open sets in the product $\prod_{\gamma\in B}(X_{\gamma}\setminus\{\xi_{\gamma}\})$, hence the sets from $\mathcal{M}'(B)$ are also open in $\prod_{\gamma\in B}X_{\gamma}$.

Let $W$ be an element of $M(B)$. Since $\mathcal{M}(B)\prec \mathcal{K}(B)$, there exists $V\in\mathcal{V}$ such that $W\subseteq V\cap F(B)$. Let \[\widetilde{W}=\pi^{-1}_{B}(\pi_{B}(W))\cap V,\] then $\widetilde{W}\cap F(B)=W$,  and $\pi_{B}(\widetilde{W})=\pi_{B}(W)$.

The set $\pi_{B}(W)$ is $\sigma$-compact as a continuous image of a  $\sigma$-compact set. The set $\pi^{-1}_{B}(\pi_{B}(W))=\pi_{B}(W)\times\sigma(X_{\gamma},\xi_{\gamma},\Gamma\setminus B)$ is $\sigma$-compact being a product of $\sigma$-compact sets. Finally, the set $\widetilde{W}$ is $\sigma$-compact as an intersection of two $\sigma$-compact sets.

Let
\begin{eqnarray*}
\mathcal{P}(B) &=& \{\widetilde{W}:W\in \mathcal{M}(B)\}; \\
\mathcal{R}(B) &=& \{\pi^{-1}_{B}(\pi_{B}(W)) : W\in \mathcal{M}(B)\}.
\end{eqnarray*}
The family $\mathcal{R}(B)$ is point-finite, because $M'(B)$ is point-finite, so $\mathcal{P}(B)$ is point-finite as a shrinking of  $\mathcal{R}(B)$. For each $W\in \mathcal{M}(B)$ we have $W\subseteq F(B)$ and $\widetilde{W}\subseteq U(B)$, therefore $\bigcup \mathcal{P}(B)\subseteq U(B)$.

Put $\mathcal{U}=\bigcup\{\mathcal{P}(B): B\in [\Gamma]^{<\omega}\}$; this family is  point-finite because the family $\{U(B):B\in [\Gamma]^{<\omega}\}$ and all families $\mathcal{P}(B), B\in [\Gamma]^{<\omega}$, are point-finite. Each set $\widetilde{W}\in \mathcal{U}$  is open in Z and $\sigma$-compact.

From the definition of sets $\widetilde{W}$, it immediately follows that $\mathcal{U}\prec\mathcal{V}$.

Finally, we shall check that $\bigcup \mathcal{U}=\bigcup\mathcal{V}$. Let $x\in V\in \mathcal{V}$. Then $x\in F(B)$, for some $B$, so $x\in V\cap F(B)\in K(B)$. Because $\mathcal{M}(B)\prec \mathcal{K}(B)$ and $\bigcup \mathcal{M}(B) = \bigcup  \mathcal{K}(B)$, there exists $W\in \mathcal{M}(B)$ such that $x\in W\subseteq V\cap F(B)$. As we observed before $\widetilde{W}\cap F(B)=W$, so $x\in\widetilde{W}$, that is $x\in \bigcup \mathcal{P}(B)$, therefore $x\in\bigcup \mathcal{U}$.
\end{proof}

Let us recall that one can adapt the Cantor-Bendixson rank of scattered spaces to
treat any topological property $\mathcal{P}$.
Namely, for any ordinal number $\alpha$ we define the $\alpha$-th  derivative
$X^{(\alpha)}$ of a given topological space $X$ with respect to  $\mathcal{P}$ as follows:

 \begin{enumerate}
\item
$X'=X^{(1)}=X\setminus\bigcup\{U\subseteq X: U \mbox{ is open and has property } \mathcal{P} \}$;
\item $X^{(\alpha+1)}=(X^{(\alpha)})'$;
\item $X^{(\alpha)}=\bigcap_{\beta<\alpha}X^{(\beta)}$ for a limit ordinal $\alpha$.
\end{enumerate}

If $X^{(\alpha)}=\emptyset$ for some $\alpha$'s then the height $ht(X)$ (with respect to $\mathcal P$)
is defined as $ht(X)=\min\{\alpha: X^{(\alpha)}=\emptyset\}$.
Note that if $X$ is compact, then $ht(X)$ is a successor ordinal number.

\begin{lem}
	\label{pochodna pozdzbioru}
	Let $\mathcal{P}$ be a hereditary topological property. For any topological spaces  $Y\subseteq X$ we have $Y^{(\alpha)}\subseteq X^{(\alpha)}$ for every $\alpha$.
\end{lem}

\begin{proof}
Note first that $Y'\sub X'$. Indeed, if $x\in Y \sm X'$ then there exists an open set
$U\subseteq X$ with property $\mathcal{P}$.
Then $U\cap Y$ also has property $\mathcal {P}$ and $x\in U\cap Y$ which implies  $x\notin Y'$.

 The general follows by induction on $\alpha$: for the successor step we apply the
 above observation while the limit case is obvious.
 \end{proof}

\begin{thm}\label{thm_d_to_a}
Let $K$ be a compact, hereditarily metacompact space and suppose that every nonempty
subspace $A$ of $K$ contains a nonempty, relatively open subspace $U$ of countable weight.
Then  $K\in\ny$.
\end{thm}

\begin{proof}
Let $\mathcal{P}$ be property of having  countable weight.
Let $X^{(\alpha)}$ be $\alpha$-th derivative of space $X$ with respect to
property $\mathcal{P}$.
Note first  that if $ht(K)=1$ then
\[K=\bigcup\{U\subseteq K: U \mbox{ is open and } w(U)\leq\omega\}.\]
By compactness,  $K$ is then  a finite union  of open sets of countable weight,
and therefore $K$ has  countable weight itself.
 Consequently, $K$ can be embedded into $I^{\omega}$, so $K\in\ny$.

Observe that, by our assumption on $K$, $ht(K)$ is well-defined so
we can check the assertion by induction on $ht(X)$.
Let $ht(X)=\alpha=\beta +1$ and assume that the theorem holds for all compact spaces of smaller height.

The space $K^{(\beta)}$ is of height 1; by the introductory remark, there is an embedding
$i:K^{(\beta)}\to I^{\omega}$. Using the Tietze-Urysohn theorem,
$i$ can be extended to  a continuous function  $f:K\to I^{\omega}$.

For every $x\in K\setminus K^{(\beta)}$ there is an open set $U_{x}\subseteq K$
such that $\overline{U_{x}}\cap K^{(\beta)}=\emptyset$.
The space $\overline{U_{x}}$ is hereditarily metacompact and,
by Lemma \ref{pochodna pozdzbioru},
$(\overline{U_{x}})^{(\beta)}\subseteq K^{(\beta)}$.
Moreover,   $\overline{U_{x}}\cap K^{(\beta)}=\emptyset$, so $(\overline{U_{x}})^{(\beta)}=\emptyset$.
It follows, that $ht(\overline{U_{x}})<\alpha$ and,   by the inductive assumption,
 $\overline{U_{x}}$ can be embedded in $\sigma(I^{\omega},\Gamma_{x})$, for some $\Gamma_{x}$.

Now  the family $\mathcal{U}=\{U_{x}:x\in K\setminus K^{\beta}\}$ forms an open cover of
$K\setminus K^{(\beta)}$. Hence, by hereditary metacompactness,
 there exists  a point-finite  open cover $\mathcal{V}=\{V_{\xi}:\xi\in \kappa\}$ of
 $K\setminus K^{\beta}$ that is inscribed in $\mathcal{U}$.
 By virtue of Lemma \ref{pokr otw}, each set $V_{\xi}$ has a point-finite cover $\mathcal{W}_{\xi}$
 consisting of open  $\sigma$-compact sets.
 The family  $\mathcal{W}=\bigcup\{\mathcal{W_{\xi}}:\xi\in\kappa\}$ is then  a point-finite cover
 of $K\setminus K^{(\beta)}$ consisting of open $\sigma$-compact sets and $\mathcal{W}$
 is  inscribed in $\mathcal{V}$.

 It follows that for every $W\in\mathcal{W}$ there is  an embedding
 \[ g_{W}:\overline{W}\to \sigma(I^{\omega},\Gamma_{W})\subseteq (I^{\omega})^{\Gamma_{W}},\]
 for some $\Gamma_W$.
  Again, by the Tietze-Urysohn theorem, $g_W$ can be extended to  a continuous function
  $G_{W}:K\to (I^{\omega})^{\Gamma_{W}}$.

Note that for any open, $\sigma$-compact set $P\subseteq K$ there exists a continuous function
$h_{P}:K\xrightarrow[]{}[0,1]$ such that $P=h_{P}^{-1}((0,1])$.
For every $W\in{\mathcal W}$ we define  $H_{W}:K\to I^{\omega}$  by
$H_{W}=h_{W} \cdot G_{W}$. Consider now the function \[T=\Big(\mathop{\BTU}_{W\in\mathcal{W}}H_{W}\Big)\bigtriangleup
\Big(\mathop{\BTU}_{W\in\mathcal{W}}h_{W}\Big)
\bigtriangleup f :K\to \sigma(I^{\omega},\Gamma).\]

Each $x\in K$ belongs to only finitely many elements of $\mathcal{W}$,
hence $H_W(x)=0=h_W(x)$ except for a finite number of $W$'s.
Consequently,
the range of  $T$ is a subset of $\sigma(I^{\omega},\Gamma)$.

To complete the proof, it remains to prove that $T$ is injective.
Consider any $x,y\in K$, $x\neq y$.

\begin{enumerate}[ {Case} (1)]
\item If $x,y\in K^{\beta}$  then $T(x)\neq T(y)$ since $f(x)\neq f(y)$.
\item If $x\in K^{(\beta)}$ and $y\notin K^{(\beta)}$ then there is
$W\in\mathcal{W}$ such that $y\in W$ and $x\notin W$;
then $h_{W}(x)=0\neq h_{W}(y)$.
\item If $x,y\notin K^{(\beta)}$ and $x, y\in W$ for some $W\in \mathcal{W}$
then either $h_{W}(x)\neq h_{W}(y)$ (and so $T(x)\neq T(y)$)
or $h_{W}(x)=h_{W}(y)$ but then  then $H_{W}(x)\neq H_{W}(y)$
because  $G_{W}(x)=g_{W}(x)\neq g_{W}(y)=G_{W}(y)$.
\item The remaining case is
$x\in W$ and $y\notin W$ for some $W\in \mathcal{W}$ (or vice versa);
 then $h_{W}(X)\neq 0=h_{W}(y)$ so again $T(X)\neq T(y)$.
\end{enumerate}
\end{proof}

\begin{proof}[Proof of Theorem \ref{charact_ec_sigma_m}]
The equivalences  (i)$\Leftrightarrow$(ii)$\Leftrightarrow$(iii)
were proved in  \cite[Theorem 1, Lemma 1]{NY}.
The implication (i)$\Rightarrow$(iv) follows from Proposition \ref{sigma-m otw zb ma otw z przel bazą} and
Lemma \ref{pokr otw}. The reverse implication is given by Theorem \refeq{thm_d_to_a}.
\end{proof}

 We can now present  the following counterpart of Theorem \ref{charact_ec_sigma_m}.

\begin{thm}\label{charact_omega_Corson}
For a compact space $K$, the following
conditions are equivalent:
\begin{enumerate}[(i)]
\item $K$ is $\omega$-Corson;
\item $K$ has a closure preserving cover consisting of finite dimensional metrizable compacta;
\item $K$ is hereditarily metacompact and each nonempty subspace $A$ of $K$ contains a nonempty relatively open, finite dimensional subspace $U$ of countable weight.
\end{enumerate}
\end{thm}

Here we again build on  the results from \cite{NY} and \cite{Ya}; in particular,
in the  proof  below we use  a slight modification of Yakovlev's reasoning from \cite[Theorem 3(a)]{Ya}.

\begin{prp}\label{prp_omega_corson_CPF}
Every $\omega$-Corson compact space $K$ has a closure-preserving cover consisting of metrizable, compact, finitely dimensional sets.
\end{prp}

\begin{proof}
By Proposition \ref{prp_sigma_prod}, we can assume that $K$ is a compact subspace of
$\sigma(I,\Gamma)$, for some $\Gamma$.
For any finite $\Gamma_{0}\sub\Gamma$ we put
\[A_{\Gamma_{0}}=\{x\in K:\forall_{\gamma\in\Gamma\setminus\Gamma_{0}}\ x(\gamma)= 0 \},\]
and consider $\mathcal{A}=\{A_{\Gamma_{0}}:\Gamma_{0}\in [\Gamma]^{<\omega)} \}$.
Then the family $\mathcal{A}$ covers $K$ and consists of finitely dimensional metrizable compacta
so it remains to prove that it is  closure-preserving.

Since   $K$ is  Eberlein compact,  it is Frechet-Urysohn so, in particular,
the closure in $K$ coincides with the sequential closure.
Therefore  it suffices to check that
the union of any subfamily  of $\mathcal{A}$ is sequentially closed.
For that purpose
consider $x_{n}\in A_{\Gamma_{n}}\in\mathcal{A}$ and assume that
the sequence of  $x_{n}$  converges to $x\in K$.

Suppose that $x\notin \bigcup_{n}A_{\Gamma_{n}}$;
then for every $n$ there is $\gamma_{n}\notin\Gamma_{n}$ such that
$x(\gamma_{n})\neq 0$. Since $x\in K$, it  has only finitely many nonzero coordinates and therefore
 there exists $\gamma\in\Gamma$ such that
the set $N=\{n:\gamma_n=\gamma\}$ is infinite.
Since $\gamma_n=\gamma\notin \Gamma_n$ for $n\in N$, we get
\[0\neq x(\gamma)=\lim_{n\in N} x_n(\gamma)=\lim_{n\in N} x_n(\gamma_n)=0,\]
a contradiction.
\end{proof}

\begin{lem}
\label{metr ośr podpokr}
Let $X$ be a separable metrizable  topological space and let $\mathcal{A}$ be a
closure-preserving cover of the space $X$ consisting of closed subsets.
Then there exists a countable subfamily $\mathcal{A}'\subseteq \mathcal{A}$ such that $\bigcup \cA'=X$.
\end{lem}

\begin{proof}
Suppose  that every countable family $\mathcal{A}'\subseteq \mathcal{A}$ is not a cover of $X$.
Then we can construct by induction on $\alpha<\omega_1$  countable subfamilies $\mathcal{A}_{\alpha}\subseteq\mathcal{A}$ such that for
$\alpha<\beta<\omega_1$ we have
$\bigcup\mathcal{A}_{\alpha}\subsetneq \bigcup\mathcal{A}_{\beta}$.
However, the space $Y=\bigcup_{\alpha<\omega_{1}}(\bigcup \mathcal{A}_{\alpha})$ is
separable and metrizable; therefore it does not allow
strictly ascending uncountable chains of closed sets, a  contradiction.
\end{proof}

\begin{proof}[Proof of Theorem \ref{charact_omega_Corson}]
The implication (i)$\Rightarrow$(ii) is given by Proposition \refeq{prp_omega_corson_CPF}.

To prove (ii)$\Rightarrow $(iii) note first that
by Theorem \ref{charact_ec_sigma_m} and Lemma \ref{pokr otw}, the
space $K$ is hereditarily metacompact.

It remains to show the second part of (iii). If
$\overline{A}$ meets the condition (iii), then so does $A$ ---
 we can assume that $A$ is closed. We can also assume that $A=K$,
because condition (ii) is inherited by closed subspaces.
By Proposition \ref{sigma-m otw zb ma otw z przel bazą}, there exists an open, nonempty set
$V\subseteq K$ such that $w(V)\leq\omega$.
The set $V$ is a Baire space as it is an open subset of a compact space.
 Let $\mathcal{A}$ be a cover as in (ii).
By virtue of Lemma \ref{metr ośr podpokr}, there exists a countable
$\mathcal{A}'\subseteq \mathcal{A}$ such that $V\subseteq \bigcup \mathcal{A}'$.
Hence, there exists $A'\in\mathcal{A}'$ such that $\wnt(A'\cap V)\neq\emptyset$ and the
set $U=\wnt(A\cap V)$ is as  required.

 For the proof of (iii)$\Rightarrow$(i) we closely  follow  the argument from the proof of
 Theorem \refeq{thm_d_to_a}. This time we consider property $\mathcal{P}$
 of having countable weight and finite dimension.

We verify (iii) by transfinite induction on  $ht(K)$ with respect to $\mathcal{P}$.
If $ht(K)=1$, then
\[K=\bigcup\{U\subseteq K: U \mbox{ is  open}, w(U)\leq\omega,\dim(U)<\infty\}.\]
Since $K$ is compact, it is a finite union of  open sets of countable weight and finite dimension;
consequently, $K$  is metrizable and finitely dimensional.
The space $K$ is $\omega$-Corson, because it can be embedded in the cube $I^{k}$,
for a natural number $k$, and of course $I^l$ can be embedded in $\sigma(I,\omega)$.

For the inductive step  one can repeat the inductive step of the proof of Theorem \refeq{thm_d_to_a},
replacing the $\sigma$-products of Hilbert cubes by the $\sigma$-products of unit intervals.
\end{proof}

\section{More properties of  $\omega$-Corson compacta and $NY$ compacta}\label{sec_more_prop}

From the characterizations of $\omega$-Corson and NY compact space given in
Section \ref{sec_charact} we can easily derive the following strengthening of
Proposition \ref{prp_metr_omega_SCD}.

\begin{cor}\label{prp_NY_omega_SCD}
A NY compact space $K$ is $\omega$-Corson if and only if it is strongly countably dimensional.
\end{cor}

\begin{proof}
The ``only if'' part follows immediately from Proposition \ref{prp_omega_SCD}.

Let $K$ be a NY compact space which is a countable union of finite dimensional compacta $A_n$.
We will show that $K$ satisfies condition (iii) of Theorem \ref{charact_omega_Corson}.
By condition (iv) of Theorem \ref{charact_ec_sigma_m},
$K$ is hereditarily metacompact.  Similarly as in the proof of the implication (ii)$\Rightarrow$(iii) of
 Theorem \ref{charact_ec_sigma_m}, we only need to verify the second part of (iii)
 for a nonempty closed subset $A$ of $K$.

By Theorem \ref{charact_ec_sigma_m}(iv), there exists a relatively open, nonempty set
$V\subseteq A$ such that $w(V)\leq\omega$. The set $V$ is a Baire space
since it is an open subset of a compact space $A$. It follows, that for some $n$, the intersection
$V\cap A_n$ has a nonempty interior in $V$.
The set $U = \wnt_V(V\cap A_n)$ being a Lindel\"of subspace of a finite dimensional compact space $A_n$,  is finite dimensional (cf.\ \cite[Theorem 3.1.23]{En}), so it has all required properties.
\end{proof}

The classes of $\omega$-Corson and NY compacta are clearly  closed under finite products;
 they are, however,  not closed under taking nontrivial countable products.
 Namely, from condition (iv) of Theorem \ref{charact_ec_sigma_m}
 we can easily deduce the following fact which also allows one  to produce
 simple examples of Eberlein compacta which are not NY compact.

\begin{cor}\label{cor_prod_not_ny}
For any sequence $(K_n)_{n\in\omega}$ of  nonmetrizable Eberlein compacta, the product
$\prod_{n\in\omega} K_n$ does not belong to $\ny$.
\end{cor}

Nakhmanson and Yakovlev noted that one can use condition (iii) of Theorem \ref{charact_ec_sigma_m}
to conclude  the following.

\begin{cor}\label{cor:ny_cont_im}
The class $\ny$ is closed under taking  continuous images.
\end{cor}

Note that the class of $\omega$-Corson compact spaces is clearly not stable under
taking Hausdorff continuous images, as the Hilbert cube is a continuous image of the Cantor set $2^\omega$.

In the context of hereditary metacompactness of NY compacta,
  the following result of Gruenhage \cite[Theorem 2.2]{Gr} is worth recalling
   (here $\Delta = \{(x,x): x\in K\}$ is the diagonal).

\begin{thm}[Gruenhage]\label{charact_Gruenhage}
	For a compact space $K$, the following
	conditions are equivalent:
	\begin{enumerate}[(i)]
	\item $K$ is Eberlein compact;
	\item $K^2$ is hereditarily $\sigma$-metacompact;
	\item $K^2\setminus\Delta$ is  $\sigma$-metacompact.
	\end{enumerate}
\end{thm}

We shall now mention  two examples showing that one cannot omit any of the two properties
named by  Theorem \ref{charact_ec_sigma_m}(iv) and Theorem \ref{charact_omega_Corson}(iii).

The first example is essentially due to Gruenhage, see his remarks
at the end of section 2 in \cite{Gr}.
He described it as an inverse limit $X= \varprojlim X_n$ of spaces $X_n$, where
$X_0 = A(\omega_1)$,  $X_{n+1}$ is obtained from $X_n$,
by replacing each isolated point of $X_n$,
by a copy of $A(\omega_1)$, and the bonding maps are the obvious quotient maps.
It is stated in \cite{Gr}, without a proof, that the space $X^2\setminus\Delta$ is metacompact.
We present below a different description of that space and prove that its every finite power  is hereditarily
metacompact (Theorem \ref{thm_prod_not_HM} shows that we cannot improve this for infinite products).
Recall here that, in general, a square $L^2$ of a hereditarily metacompact (even hereditarily Lindel\"of)
compact space $L$ need not to be hereditarily metacompact, which is witnessed by
the \emph{double arrow} space $L$ (cf.\ \cite[53.B(a)]{En1}).

Recall also that \emph{polyadic} spaces, introduced by Mr\'owka, are the continuous images of  spaces of the form $A(\kappa)^\lambda$.

\begin{ex}\label{ex_Grunehage}
There exists a zero-dimensional  uniform Eberlein and polyadic compact space $K$
such that  $K\notin\ny$ but  $K^n$ is hereditarily metacompact for every $n$.
\end{ex}

\begin{proof}
Consider the tree $T = \omega_{1}^\omega \cup \bigcup_{n\in\omega} \omega_{1}^n$  with the standard order $\preceq$ given by inclusion. For $s\in \omega_{1}^n,\ n\in\omega$, let $V_s = \{t\in T: s\preceq t\}$ (a wedge in $T$).

Our example $K$ is the tree $T$ equipped with the coarse wedge topology, i.e., the topology generated by the sets $V_s,\ T\setminus V_s$ for $s\in \omega_{1}^n,\ n\in\omega$. Since $T$ has only one minimal element, and every branch in $T$ has a greatest element, the space $K$ is compact (cf. \cite[Theorem 3.4]{Ny2}). From the definition of the topology on $T$ it immediately follows that $K$ is zero-dimensional.

To show that $K$ is uniform Eberlein compact and polyadic it is enough to check that it is a continuous image of the product $A(\omega_{1})^\omega$. Indeed, one can easily verify that the map $\varphi: A(\omega_{1})^\omega\to K$ defined by
\[\varphi((x_n)_{n=0}^\infty) =
\begin{cases}
(x_i)_{i=0}^\infty& \mbox{if } x_i \in\omega_1 \mbox{ for all } i\in\omega\\
(x_i)_{i=0}^{n-1}& \mbox{if } x_n =\infty \mbox{ and } x_i \in\omega_1 \mbox{ for all } i < n\,,
\end{cases}\]
for $(x_i)_{i=0}^\infty \in A(\omega_{1})^\omega$, is a continuous surjection.

Let
\[\mathcal{B} = \{V_s\setminus \bigcup_{i=0}^k V_{t_i}: s\in \omega_{1}^n,\ t_0,\dots,t_k\in  \omega_{1}^{n+1},\ n,k\in\omega\}\,.\]
The family $\mathcal{B}$ is a base of $K$ having the following properties
(observe that (1)$\Rightarrow$(2) and (3)$\Rightarrow$(4)):

\begin{enumerate}
\item $u \in V_s\setminus \bigcup_{i=0}^k V_{t_i}$ if and only if $s\preceq u$ and $t_i\npreceq u$ for $i\le k$;
\item for a fixed $u\in T$, the collection $\mathcal{B}_u = \{U\in \mathcal{B}: u\in U\}$ is closed under finite unions;
\item $V_s\setminus \bigcup_{i=0}^k V_{t_i} \subseteq V_{s'}\setminus \bigcup_{i=0}^{k'} V_{t'_i}$
 implies $s'\prec s$ or $s' = s$ and $\{t'_0,\dots,t'_{k'}\} \subset \{t_0,\dots,t_k\}$;
\item for a fixed $u\in T$, the collection $\mathcal{B}_u$ does not contain infinite strictly increasing chains.
\end{enumerate}

Properties (2) and (4) imply that the base $\mathcal{B}$ is
\emph{point-additively Noetherian} (cf.\ \cite[Definition 1]{Ny1}).
Therefore, Theorems 1 and 5 from \cite{Ny1} yield that every finite product $K^n$ is hereditarily metacompact.

Finally, the space $K$ does not belong to $\ny$, since each nonempty open subset of $K$ contains a copy of $A(\omega_{1})$, which violates the conclusion of Proposition \ref{sigma-m otw zb ma otw z przel bazą}.
\end{proof}

\begin{ex}\label{ex_non_Eberlein}
For the second example consider any scattered compact space $K$  which is not Eberlein compact.
Then $K$ has the property that each nonempty subspace $A$ of $K$ contains a nonempty
relatively open separable, metrizable, finite dimensional subspace $U$.
\end{ex}

Another relevant example comes from \cite{Ya}:

\begin{ex}\label{ex_Yakovlev}
The product $A(\omega_1)^\omega$ is Eberlein compact, but not hereditarily metacompact.
\end{ex}

As we show below, such a property is shared by all infinite product of nonmetrizable Eberlein compacta,
see Theorem \refeq{thm_prod_not_HM}.

\begin{lem}
\label{Wn z Worrella}
Let $K$ be a compact, hereditarily metacompact topological space. If
$\vf :K \to L$ is a continuous surjection  then $L$ is hereditarily metacompact as well.
\end{lem}

\begin{proof}
For any $B\subseteq L$,
by our assumption, the space $A= \phi^{-1}[B]$ is metacompact.
By compactness of $K$,  the mapping $\vf$ is closed and so is
the restriction of $\vf$ to $A$
 (if $F$ is closed in $A$ then $F=A\cap H$ for some closed subset $H$ of $K$;
 then $\vf[F]=\vf[H]\cap B)$).
By Worrell's theorem, see \cite[Theorem 5.3.7]{En1},
 metacompactness is preserved by closed mappings
so $B$ is metacompact.
\end{proof}

\begin{thm}\label{thm_prod_not_HM}
If $(K_{n})_{n\in\omega}$ is a sequence of nonmetrizable Eberlein compact spaces
then the product $\prod_{n\in\omega}K_{n}$ is not hereditarily metacompact.
\end{thm}

\begin{proof}
It is well-known  that a nonmetrizable Eberlein compactum is not $ccc$, so for each $n$,
 there exists a family $\{U_n^{\alpha}:\alpha\in\omega_{1}\}$ of pairwise disjoint, nonempty, open sets of $K_{n}$. Pick any $x^{\alpha}_n\in U_n^{\alpha}$. Consider
\[L_{n}=(K_{n}\setminus \bigcup_{\alpha\in\omega_{1}}U_n^{\alpha})\cup\{x_n^{\alpha}:\alpha\in\omega_{1}\}.\]
One can easily verify that the quotient space $L_{n}/(K_{n}\setminus \bigcup_{\alpha<\omega_{1}}U_n^{\alpha})$ obtained by identifying all points in $K_{n}\setminus \bigcup_{\alpha<\omega_{1}}U_n^{\alpha}$ is homeomorphic to $A(\omega_{1})$. Pick a homeomorphism $h_n$ witnessing this. Let $q_{n}: L_{n}\to L_{n}/(K_{n}\setminus \bigcup_{\alpha<\omega_{1}}U_n^{\alpha})$ be the quotient map, and let $f_n = h_n\circ q_n$. The map \[\prod_{n\in\omega}f_{n}:\prod_{n\in\omega}L_{n}\to A(\omega_{1})^{\omega}\]
is a continuous surjection, and therefore, by Example \ref{ex_Yakovlev} and  Lemma \ref{Wn z Worrella},
 the product $\prod_{n\in\omega}K_{n}$ is not hereditarily metacompact.
\end{proof}

Recall that the preorder $\le^*$ on $\omega^\omega$ is defined by
$f\le^* g$ if $f(n)\le g(n)$ for all but finitely many $n\in\omega$. A subset $A$ of $\omega^\omega$ is called \emph{unbounded} if it is unbounded with respect to this preorder.
We  refer  below to the classical cardinal number:
\begin{eqnarray*}
	\mathfrak{b} = \min\{|A|: A \text{ is an unbounded subset of } \omega^\omega\}\,.
\end{eqnarray*}
It is well known that, for every  $n\ge 1$, the statement $\mathfrak{b} = \omega_n$ is consistent with \textsf{ZFC}, (cf.\ \cite[Theorem 5.1]{vD}).

The next result was proved in \cite{Ma} for a proper subclass $\eco$ (see Section \ref{sec_eco}) of the class of NY compacta. We do not know if it can be proved for all Eberlein compact spaces.

\begin{thm}\label{thm_on_ec_sigma_m}
Assuming that $\mathfrak{b} > \omega_1$, each nonmetrizable compact space $K\in\ny$ contains a closed nonmetrizable zero-dimensional subspace $L$.
\end{thm}

\begin{proof}
Our argument follows closely the proof of Theorem 4.15 in \cite{Ma}.

Let $K\in\ny$; we can assume that $K\sub \sigma(I^\omega,\Gamma)$
for some set $\Gamma$. Since $K$ is nonmetrizable, obviously the set $\Gamma$ must be uncountable.
We can also assume that, for each $\gamma\in \Gamma$, there is $x_\gamma\in K$ such that $x_\gamma(\gamma) \ne (0,0,\dots)$. Given $\gamma\in \Gamma$, the set
\[F_\gamma = \{\delta\in \Gamma: x_\gamma(\delta) \ne (0,0,\dots)\}\]
is finite and nonempty. Using the $\Delta$-system lemma we can find a finite set $A\subseteq \Gamma$ and a set $S\subseteq \Gamma$ of size $\omega_1$ such that, for any distinct $\alpha,\beta \in S$, $F_\alpha\cap F_\beta = A$.  Now, we can identify the product $(I^\omega)^\Gamma$ with the product $I^{\omega\times\Gamma}$, and apply Lemma 4.10 from \cite{Ma} for the sets $X=\{x_\gamma: \gamma\in S\}$ and $\Gamma_{0} = \omega\times A$.
\end{proof}

\begin{rem}\label{conv_no_NY}
It follows from  Proposition \ref{sigma-m otw zb ma otw z przel bazą} 
that no nonmetrizable, compact convex subset $K$ of a topological vector space $E$ is NY compact because
every nonempty open set $V\subseteq K$ contains a copy of the space $K$. 
Indeed, for a fixed $x_0\in V$, by the continuity of the linear operations in $E$, and the compactness of $K$, the set $\{(1-\eps)x_0 +\eps y: y\in K\}$ - an affine copy of $K$, is contained in $V$, for suitably small $\eps >0$.

In particular, the space $P(K)$ of all probability Radon measures on a compact space $K$ is NY compact if and only if $K$ is metrizable (recall that $P(K)$ is Eberlein compact if and only if $K$ is Eberlein compact).
\end{rem}

\section{Between  $\omega$-Corson and $NY$ compacta}
\label{sec_eco}

One of the motivation for this section is to correct two false comments from the first author's paper \cite{Ma}.

We say that a compact space $K$ belongs to the class $\eco$ if, for some set $\Gamma$ there is an embedding $\phe:K\to \mathbb{R}^\Gamma$ and a countable subset $\Gamma_0$ of $\Gamma$ such that, for each $x\in K$, the set $\supp(\phe(x))\setminus \Gamma_0$ is finite (cf. \cite{Ma}).

\begin{prp}
A compact space $K$ belongs to $\eco$ if and only if it can be embedded into the product $I^\omega \times L$ of the Hilbert cube $I^\omega$ and some $\omega$-Corson compact space  $L$.
\end{prp}

\begin{prp}
Each $\omega$-Corson compact space belongs to the class $\eco$, and each member of $\eco$ is NY compact.
\end{prp}

Again, the Hilbert cube witnesses the fact that the first implication in the above proposition cannot be reversed. The same is true for the second implication, see Example \ref{ex_cont_im_eco} below and Proposition \ref{prp_discrete_unions}.
\medskip

In \cite{Ma} we made, without any justification, the following comment ``One can even show that, for metrizable compacta $M_n$, the countable product of Aleksandrov duplicates $AD(M_n)$ is in $\eco$.'' Corollary \ref{cor_prod_not_ny} shows that this does not hold true. Indeed,
the product $[AD(2^\omega)]^\omega$ is an example of a first-countable uniform Eberlein compact space which does not belong to the class $\ny$, hence is not in $\eco$ (cf.\ \cite[Example 4.16]{Ma}).
\medskip

The example below shows that another unexplained comment from \cite{Ma}, claiming the any member of the class $\eco$, is uniform Eberlein is also false.

\begin{ex}\label{ex_non_uniform}
Let $K$ be a scattered Eberlein compact space which is not uniform Eberlein compact. Then $K$ is $\omega$-Corson (hence belongs to $\eco$).
\end{ex}

As we mentioned in the previous section, the class of $\omega$-Corson compact spaces is not stable under taking Hausdorff continuous images. One can prove the class $\eco$ also lacks this property.

\begin{ex}\label{ex_cont_im_eco}
The space $\alpha(D(\omega_1)\times 2^\omega)$ is $\omega$-Corson (hence belongs to $\eco$), but its continuous image $\alpha(D(\omega_1)\times I^\omega)$ does not  belong to $\eco$.
\end{ex}

\begin{proof}
The fact that  $\alpha(D(\omega_1)\times 2^\omega)$ is $\omega$-Corson is given by Proposition \refeq{prp_discrete_unions}.

Suppose that there is an embedding $\phe: \alpha(D(\omega_1)\times I^\omega) \to \mathbb{R}^\Gamma$ and a countable subset $\Gamma_0$ of $\Gamma$ such that, for each $x\in \alpha(D(\omega_1)\times I^\omega)$, the set $\supp(\phe(x))\setminus \Gamma_0$ is finite. Let
\[A = \{x\in \alpha(D(\omega_1)\times I^\omega): \phe(x)|\Gamma_0 = \phe(\infty)|\Gamma_0 \}\,,\]
where $\infty$ denotes the point at infinity of  $\alpha(D(\omega_1)\times I^\omega)$. Since $\Gamma_0$ is countable, the set $A$ is a $G_\delta$-subset of $\alpha(D(\omega_1)\times I^\omega)$ containing $\infty$, hence it contains a copy of the Hilbert cube $I^\omega$. On the other hand, the projection $\pi: \mathbb{R}^\Gamma \to \mathbb{R}^{\Gamma\setminus \Gamma_{0}}$ restricted to $\phe(A)$ is an embedding and its image is contained in the $\sigma$-product of intervals, which contradicts the fact that $I^\omega$ is not strongly countably dimensional.
\end{proof}

We conclude the section with an observation related to the wider class of $\kappa$-Valdivia compact spaces
considered by Kalenda \cite{Ka}. The following provides a positive
answer to a question posed by Ond\v{r}ej Kalenda in a conversation.

\begin{prp}
For every metrizable compactum $K$ there is an embedding $f :K\to I^\omega$
such that $f(K)\cap \sigma(I,\omega)$ is dense in $f(K)$.
\end{prp}

\begin{proof}
Let $K$ be a metrizable compactum. We can assume that $K\subseteq I^\omega$.
Let $D$ be a countable, dense subset of $K$, and $E$ be a countable, dense subset of $I^\omega$, contained in $\sigma(I,\omega)$.
Put $F = D\cup E$. The Hilbert cube $I^\omega$ is \emph{countable dense
homogeneous}, i.e.,
given countable dense subsets $A,A' \subseteq I^\omega$ , there is a homeomorphism $h: I^\omega\to
I^\omega$ such that $h(A) = A'$. Let $h$ be such homeomorphism applied for $A = F$ and $A' = E$. Then the restriction $h|_K$ is an embedding such that $h|_K(K)\cap \sigma(I,\omega)$ contains $h(D)$ hence is dense in $h|_K(K)$.
\end{proof}

\section{Some local properties}\label{local_properties}

We consider here  $\kappa$-Corson compacta for $\kappa>\omega$ and discuss
their properties related to the usual cardinal invariants of topology.

Given a (topological) space $K$ and  $x\in K$, $\chi(x,K)$ denotes the local character
of the space $K$ at $x$ (i.e. the minimal cardinality of a local base). Recall that
for every point $x$ in a compact space $K$, $\chi(x,K)$ is just the minimal size of a family of open
sets which intersection is $\{x\}$.
We write $\pi\chi(x,K)$ for the $\pi$-character
of the space $K$ at $x$ (i.e. the minimal cardinality of a local $\pi$--base).
Similarily,  $t(x,K)$ denotes
the tightness of $K$ at $x$ (recall that $t(x,K)\le\tau$ means that whenever $X\sub K$ and
$x\in \overline{X}$, then there is $A\sub X$ such that $|A|\le\tau$ and $x\in\overline{A}$).
The tightness of a  space is $\le\tau$ if
$t(x,K)\le\tau$ for every $x\in K$. The strong tightness of $K$ is $\le\tau$ if
for every  $X\sub K$ and
$x\in \overline{X}$  there is $A\sub X$ such that $|A|<\tau$ and $x\in\overline{A}$.

Recall first that the following general fact holds even for singular $\kappa$,
compare \cite[Proposition 2.19]{BKT}.

\begin{lem}\label{character}
Given an uncountable $\kappa$, every $\kappa$-Corson compact space $K$ contains
a point of character $<\kappa$.
\end{lem}

\begin{proof}
Assume that  $K\sub\sikg$ and consider the pointwise order $\leq$ on elements of $K$.
It is easy to check that there is $x\in K$ which is pointwise maximal, that is $y\in K$ and
$y(\gamma)\ge x(\gamma)$ for every $\gamma\in\Gamma$ imply $y=x$.
Then $x$ is the only point in the intersection of the family of open sets
 \[ \{y\in K: y(\gamma)>x(\gamma)-1/n\}, \gamma\in\supp(x), n\in\en,\]
  so there is a local base at $x\in K$ of size $<\kappa$.
\end{proof}

We can now  answer Question 1.22(ii) in \cite{Ka}; the fact given below was also
noted in \cite[Corollary 2.20]{BKT}.

\begin{cor}
The space $[0,1]^{\kappa}$ is $\kappa$--Corson for no $\kappa\ge \omega$.
\end{cor}

\begin{proof}
For $\kappa=\omega$ this follows by Proposition \ref{prp_omega_SCD}.
For $\kappa>\omega$ we can use Lemma \ref{character}, since
$\chi(x,[0,1]^{\kappa})=\kappa$ for every $x\in [0,1]^{\kappa}$.
\end{proof}

To discuss the tightness we first recall  \cite[Lemma 1.19]{Ka}.

 \begin{lem} \label{st}
 For every $\kappa>\omega$ and any $\Gamma$ the space $\sikg$ has strong tightness $\kappa$ if and only if $\kappa$ is regular.
  \end{lem}

\begin{lem}\label{picharacter}
If $\kappa>\omega$ is regular, then $\pi\chi(x,K)<\kappa$ for every $\kappa$-Corson compact space
$K$ and every $x\in K$.
\end{lem}

\begin{proof}
Note first that Lemma \ref{character} in fact implies that the points $x\in K$ satisfying
$\chi(x,K)<\kappa$ form a dense subset $D$ of  $K$.
By Lemma \ref{st}, for every $y\in K$ there is $D_0\sub D$ such that $|D_0|<\kappa$ and $y\in\overline{D_0}$.
For every $x\in D_0$ pick a local base $\cU(x)$ at $x$ of size $<\kappa$.
Then $\bigcup_{x\in D_0} \cU(x)$ is a $\pi$-base at $y$ of size $<\kappa$.
\end{proof}

The argument above is not applicable in the case of singular $\kappa$; nonetheless, the
same assertion holds.
The following result answers Question 1.22(1) from \cite{Ka}; here we
use an idea due to   Tka\v{c}enko, see  \cite[Lemma 1]{Tk83}.

\begin{thm}\label{tkacenko}
  Let $\kappa$ be an uncountable singular
cardinal, and let $K$ be $\kappa$--Corson compact.
Then $\pi\chi(x,K)<\kappa$ for every $x\in K$.

Consequently,  $K$ has strong tightness $\le\kappa$.
\end{thm}

\begin{proof}
Consider again a compact space $K\sub\sikg$.
Let $\lambda={\rm cf}(\kappa)$, and let $\vf:\lambda\to\kappa$ be an increasing
function with an unbounded range. For any $\xi<\lambda$ we put
\[A_{\xi}=\{x\in K:\; \pi\chi(x,K)\le |\vf(\xi)|\mbox{ and } |\supp(x)|\le |\vf(\xi)|\}.\]

\noindent {\sc Claim 1.}
  For every $\xi<\lambda$ and $x\in K$, if $x\in\overline{A_{\xi}}$, then $\pi\chi(x,K)<\kappa$.
\medskip

Indeed, let $\tau=\max(|\vf(\xi)|,|{\rm supp}(x)|)$. By Lemma \ref{st},  the space
$\Sigma_{\tau^{+}}(\er^{\Gamma})$ has tightness $\tau$.
Hence $x\in\overline{A_{\xi}}$ lies
 in the closure
of $\tau$ many points $y$ with $\pi\chi(y,K)\le\tau$. This clearly implies that
 $\pi\chi(x,K)\le\tau<\kappa$.
\medskip

\noindent {\sc Claim 2.}
$K=\bigcup_{\xi<\lambda}\overline{A_{\xi}}$.
\medskip

Suppose otherwise; let $K'=\bigcup_{\xi<\lambda}\overline{A_{\xi}}$ and $x\in K\sm K'$.
For every $\xi<\lambda$ we may find a closed $G_\delta$ set $Z_{\xi}\sub K$ such that $x\in Z_{\xi}$
and $Z_{\xi}\cap \overline{A_{\xi}}=\emptyset$. Write $Z=\bigcap_{\xi<\lambda}Z_{\xi}$.
Then $Z$ is nonempty and, by Lemma \ref{character}, there is $z\in Z$ such that $\chi(z,Z)<\kappa$.
Choose  $\xi$ such that
\[ \vf(\xi)\ge \max (\chi(z,Z), |\supp(z)|,\lambda) .\]
 Note that
\[ \pi\chi(z,K)\le \chi(z,K)\le \chi(z,Z)\cdot \lambda,\]
where the latter inequality follows from the fact that $Z$ is an intersection of $\lambda$ many open sets.
But then $z\in A_\xi$,
a contradiction with $z\notin K'$.
\medskip

It follows immediately from Claim 1 and Claim 2 that $\pi\chi(x,K)<\kappa$ for every $x\in K$.
This clearly implies the final statement, for
if $L=\overline{X}\sub K$ and $ x\in L$, then $L$ is $\kappa$--Corson, so
$\pi\chi(x,L)<\kappa$, which implies  that there is $Y\sub X$ with $|Y|<\kappa$ and $x\in \overline{Y}$.
\end{proof}

As we shall see later, it seems to be unclear if the class of $\kappa$-Corson
compacta is stable under continuous images in case of singular $\kappa$.
For that reason we note the following.

\begin{cor}\label{witek}
For every uncountable $\kappa$, if $L$ is a continuous image of a $\kappa$-Corson
compact space, then $\pi\chi(y,L)<\kappa$ for every $y\in L$.
\end{cor}

\begin{proof}
Let $f:K\to L$ be a continuous surjection where $K$  $\kappa$-Corson compact.
It is well-known that there is closed $K_0\sub K$ such that
$f[K_0]=l$ and $f|K_0$  is irreducible,  i.e.\ $f[U]$ has a nonempty interior
for every nonempty open set $U\sub K_0$. Then $K_0$ is again $\kappa$-Corson so
it is enough to consider the case $K_0=K$.

By Lemma \ref{picharacter}
and Theorem \ref{tkacenko}, $\pi\chi(x,K)<\kappa$ for every $x\in K$.
It follows easily from irreducibility that $\pi\chi(f(x),L)\le\pi\chi(x,K)<\kappa$.
\end{proof}

The following simple example comes from \cite[1.21]{Ka} and exhibits another
difference between regular and singular $\kappa$.

\begin{ex}\label{ex1}
Fix an uncountable $\kappa$, and let
 $[0,\kappa]$ denote the space of ordinal numbers $\le\kappa$ equipped with the
order topology. The space $ [0,\kappa]$ is $\kappa$-Corson
if and only if  $\kappa$ is singular.
\end{ex}

Indeed, if $\kappa$ is regular, then examining $\kappa\in\overline{[0,\kappa)}$
we see that the space does not have strong tightness $\kappa$.
On the other hand, if $\kappa$ is singular, then take a cofinal set $C\sub\kappa$ with $|C|<\kappa$.
Now the family $\mathcal U$ of clopen sets $(\alpha,\beta]$, where $\beta\in C$ and $\alpha<\beta$
separates points of $[0,\kappa]$ and $\ord(x,\mathcal{U}) < \kappa$ for all $x\in [0,\kappa]$
so $[0,\kappa]$ is $\kappa$-Corson by
Proposition \ref{charact_BKT}.

We enclose here some remarks on converging sequences.

\begin{lem}\label{cs}
If a $\kappa$-Corson  compact space  $K$ contains a point $x_0$  such that
 $\chi(x_0,K)=\kappa$, then $K$ contains a nontrivial converging sequence.
\end{lem}

\begin{proof}
Consider a compact space  $K\sub\Sigma_\kappa([0,1]^\Gamma).$

Write   $S_0=\supp(x_0)$;
as $\chi(x_0,K)=\kappa$, there is
$x_1\in K$ such that $x_1|S_0=x_0|S_0$, but $x_1\neq x_0$. Then we consider the
set $S_1=S_0\cup \supp(x_1)$ and repeat the argument.

In this fashion, we can define a sequence of  points $(x_n)_{n\ge 1}$ in $K$,
and  $S_0\sub S_1\sub\ldots\sub\Gamma$ such that for every $n\ge 1$

\begin{itemize}
\item[--] $x_n\neq x_0$,
\item[--] $x_n|S_n\equiv x_0|S_n$, and
\item[--] $\supp(x_n)\sub S_n$.
\end{itemize}

Then $x_n$ converge to $x_0$: if $\gamma\in\bigcup_n S_n$, then $x_n(\gamma)=x_0(\gamma)$
eventually holds; otherwise $x_n(\gamma)=0=x_0(\gamma)$.
\end{proof}

Note, however, that if the character of every point in a $\kappa$-Corson space is $<\kappa$, then
$K$ may contain
no converging sequence.
Such a phenomenon occurs for $\kappa=\cont^+$ when we let $K$ be $\beta\omega$ (see Remark \ref{basic_observations}(c)).

\section{Good sets}\label{good_sets}
The concept of a good set of coordinates was widely used in analysing
the properties of  classical Corson compacta, see \cite{Ne}, \cite{AMN88} and \cite{Ka}.
The same concept can be applied to study $\kappa$-Corson compacta
for a regular uncountable $\kappa$.

Consider any fixed $K\sub [0,1]^\Gamma$; a set $A\sub K$ is determined by coordinates
in $H\sub \Gamma$ if for every $x\in A$ and $y\in K$,  $x|H=y|H$ implies $y\in A$; we write
$A\sim H$ in such a case.
In the sequel, by a basic open set in $K$  we mean a set of the form
\[ \{y\in K: |y(\gamma)-\phi(\gamma)|<1/n\mbox{ for all } \gamma\in I\},\]
where $I\sub\Gamma$ is finite and $\phi:I\to\qu$.
The following fact is standard.

\begin{lem}\label{l4}
If $K$ is a compact subset of  $[0,1]^\Gamma$, then
every clopen subset of $K$ is determined by finitely many coordinates.
\end{lem}

If $K\sub \Sigma_\kappa([0,1]^\Gamma)$ and $H\sub\Gamma$, then we write
\[ \pi_H: K\to [0,1]^\Gamma,\]
for the corresponding `projection'; precisely $\pi_H(x)(\gamma)=x(\gamma)$ if $\gamma\in H$ and
$\pi_H(x)(\gamma)=0$ otherwise.

\begin{definition}
Given $K\sub \Sigma_\kappa([0,1]^\Gamma)$,
a set $H\sub \Gamma$ is said to be {\em good for} $K$ if $\pi_H[K]\sub K$.
\end{definition}

The following observation follows easily by compactness.

\begin{lem}\label{la}
Given a compact space $K\sub \Sigma_\kappa([0,1]^\Gamma)$, the union of any chain of subsets of $\Gamma$
which are good for $K$ is also good.
\end{lem}

The proof of the next lemma closely follows Benyamini's argument for
$\kappa=\omega_1$, see
\cite[6.24]{Ne}.

\begin{lem}\label{l2a}
Suppose that  $\kappa$ is a regular uncountable cardinal and let $\lambda\ge\kappa$.
Given a compact space $K\sub \Sigma_\kappa([0,1]^\lambda)$,
 for every
$H\sub\lambda$ with $|H|<\lambda$   there is $G(H)\supseteq H$ such that
$|G(H)|<\lambda$ and $G(H)$ is good for $K$.
\end{lem}

\begin{proof}
Let $H_0=H$; we construct and increasing sequence  $H_n\sub  \lambda$ as follows.
Given $H_n$, for every nonempty basic open set $U$ in $K$ depending on coordinates in $H_n$
we pick $x_U\in U$ and set
\[H_{n+1}=H_n\cup \bigcup_{U\sim H_n} \supp(x_U).\]
Then we put $G(H)=\bigcup_n H_n$.

If $\lambda=\kappa$, then $|H_n|<\kappa$ and $|G(H)|<\kappa$ by regularity of $\kappa$.
If $\lambda>\kappa$, then writing $\tau=\max(\kappa, |H|)$, we have
$|H_n|\le\tau$ and so $|G(H)|\le\tau<\lambda$.

Take any $x\in K$ and consider $z=\pi_{G(H)}(x)$. If $V$ is a basic neighborhood of
$z\in\sikl$ and $V\sim J\sub\lambda$, then for $I=J\cap G(H)$ we have $I\sub H_n$ for some $n$
and we can write $V=U_1\cap U_2$, where $U_1\sim I$, $U_2\sim J\sm I$.
Hence there is $x\in K$ with $\supp(x)\sub H_{n+1}\sub G(H)$ such that
$x\in U_1\cap U_2=V$. This shows that $z\in \overline{K}=K$, and we are done.
\end{proof}

\begin{lem}\label{l2b}
Suppose that  $\kappa$ is a regular uncountable cardinal and let $\lambda\ge\kappa$.
Given a compact space $K\sub \Sigma_\kappa([0,1]^\lambda)$,
and a continuous surjection $f:K\to L$,
 for every
$H\sub\lambda$ with $|H|<\lambda$   there is $G^f(H)\supseteq H$ such that
\begin{enumerate}[(i)]
\item $|G^f(H)|<\lambda$ and $G^f(H)$ is good for $K$.
\item for every $x,y\in K$, if $f(x)=f(y)$, then $f(\pi_{G^f(H)}(x))=f(\pi_{G^f(H)}(y))$.
\end{enumerate}
\end{lem}

\begin{proof}
Firstly, consider a good set $G\sub\lambda$
and let $K_0=\pi_G[K]\sub K$.
Suppose that $U_1, U_2$ are basic open sets in $K_0$ depending on $G$
and $f[\overline{U_1}]\cap f[\overline{U_2}]=\emptyset$.
\medskip

\noindent{\sc Claim.}
There is a finite set $I=I(U_1,U_2)\sub\lambda$
and  open sets $V_1,V_2$ in $K$ determined by $I$ such that
\[ V_i\supseteq \overline{U_i}\mbox{ for } i=1,2 \mbox{ and } f[V_1]\cap f[V_2]=\emptyset.\]

Indeed, take disjoint open sets $W_i$ in $L,\ i=1,2$, separating $f[\overline{U_1}]$ and $f[\overline{U_2}]$, and find sets  $V_i$ which are finite unions of basic open sets determined by a finite set $I\subseteq \lambda$, and satisfy
$\overline{U_i} \subseteq V_i \subset f^{-1}(W_i)$ for $i=1,2$.
\medskip

Using Claim we define $G_0\sub G_1\sub\ldots \lambda$ as follows:
Start  with $G_0=G(H)$; given $G_n$ we let  $H_n$ be  the union of $G_n$
and all the finite sets $I(U_1,U_2)$ where $U_1,U_2$ are as in Claim.
Then we put $G_{n+1}=G(H_n)$.

Finally, $G^f(H)=\bigcup_n G_n$ is as required. Indeed, such a set is good for $K$
by Lemma \ref{la}. Suppose $x,y\in K$ are such that $f(x')\neq f(y')$, where
$x'=\pi_{G^f(H)} (x)$, $y'=\pi_{G^f(H)} (y)$. Then there are basic neighborhoods
$U_1,U_2$ of $x',y'\in \pi_{G^f(H)}[K]$ such that $f[\overline{U_1}]\cap f[\overline{U_2}]=\emptyset$.
Then $U_1,U_2$ are determined by a finite set $J\sub G^f[H]$;
we have $J\sub G_n$ for some $n$ and our  construction at this step added coordinates showing that
 $f(x)\neq f(y)$.
\end{proof}

For the sake of section \ref{bds} we need to rephrase Lemma \ref{l2a} as follows
(here the proof is very similar).

\begin{lem}\label{l2}
Suppose that  $\kappa$ is a regular cardinal.
Given a compact space $K\sub \Sigma_\kappa([0,1]^\kappa)$,
 for every
$H\sub\kappa$ of size $<\kappa$ there is $\xi<\kappa$ such that
$H$ is contained in the interval $ \{\beta:\beta<\xi\}$ which is good for $K$.
\end{lem}

In section \ref{bds} we also use the following easy observation.

\begin{lem}\label{l4.5}
Let $K$ be a compact subset of  $\Sigma_\kappa([0,1]^\Gamma$) and let $H\sub\Gamma$ be a good set for $K$;
write $K'=\pi_H[K]\sub K$. If $U$ is  a clopen subset of $K$, then $U\cap  K'=C\cap K'$ for some clopen set $C$ determined by a finite subset of $H$.
\end{lem}

\section{Continuous images}\label{cont_images}

We  have already seen that  the class $\ny$ is stable under continuous images while
$\omega$-Corson compacta do not have such a property, see Corollary \ref{cor:ny_cont_im} and
the remark after it. Then the following theorem resolves the case $\kappa=\omega_1$.

\begin{thm}\label{ci:1}
The class of Corson compact spaces is stable under continuous images.
\end{thm}

The result  was first  proved by Gul'ko \cite{Gu77} and
by Michael \& Rudin \cite{MR77};
the proof in \cite{MR77} is a modification of  the authors' new argument for stability of the class
of Eberlein compacta under continuous images. However, according to
Kalenda \cite[page 2]{Ka} that modification is not quite straightforward.
Quite different proofs of \ref{ci:1} were given by Gruenhage \cite{Gr} (using
covering properties and topological games) and by Pol \cite{Po84}
(using a characterization of Corson compacta $K$ via properties of
$C_p(K)$; see \cite{Ka99} or \cite{Ka}).

Bell and Marciszewski \cite[Theorem 7.1]{BM04} extended Pol's approach
and proved that if $K$ is $\tau^+$-Corson, then so is  every
continuous image of $K$.

The following is Corollary 2.17 in \cite{BKT}.

\begin{thm}\label{stable}
For every regular cardinal number $\kappa$,
the class of $\kappa$-Corson compacta is stable under continuous images.
\end{thm}

That result was obtained via model-theoretic approach. We give below a purely
topological argument for \ref{stable} --- it closely follows
basic ideas from Gul'ko \cite{Gu77} (mentioned also in \cite[6.26]{Ne}).
Namely, we shall prove by induction on $\lambda$ the following.

\begin{thm}\label{stable_ind}
Given a regular cardinal number  $\kappa$,
for every $\lambda$, if $K\sub\sikl$ is compact and $L$ is a continuous image of $K$,
then $L$ is $\kappa$-Corson compact.
\end{thm}

\begin{proof}
Fix  a compact space $K\sub \sikl$ and a continuous surjection $f:K\to L$.

By Lemma \ref{l2b} (and Lemma \ref{la}),  there is a family
$\{G_\xi:\xi<\lambda_0 \}$ (where $\lambda_0\le\lambda$,
 note that  $\lambda$ may be singular)
of nonempty and strictly increasing subsets of $\lambda$ such that $\lambda=\bigcup_{\xi<\lambda} G_\xi$ and for every $\xi<\lambda_0$

\begin{enumerate}[(i)]
\item $|G_\xi|<\lambda$;
\item $G_\xi$ is good for $K$;
\item for every $x,y\in K$, if $f(x)=f(y)$, then $f(\pi_{G_\xi}(x))=f(\pi_{G_\xi}(y))$;
\item $G_\xi=\bigcup_{\eta<\xi} G_\eta$ whenever $\xi$ is a limit ordinal.
\end{enumerate}

Write $p_\xi= \pi_{G_\xi}$, $K_\xi=p_\xi[K]$ and $p^\eta_\xi$ for $p_\xi$ restricted to
$K_\eta$ whenever $\xi<\eta$.
Note that every $p_\xi$ is a retraction and $p_\xi\circ p_\eta=p_\xi$ for $\xi<\eta$, so
we obtain the following inverse system

\begin{center}
\begin{tikzcd}
K_0 & \arrow[l, "p^\xi_0" above]  K_\xi & \arrow[l, "p^\eta_\xi" above]   K_\eta \arrow[l]  & \arrow[l, "p_\eta" above] K
\end{tikzcd}
\end{center}

Using $(ii)$  above we will get the following commutative diagram

\begin{center}
\begin{tikzcd}
K_0\arrow[d,"f_0"] & \arrow[l, "p^\xi_0" above]  K_\xi\arrow[d, "f_\xi"] & \arrow[l, "p^\eta_\xi" above] \arrow[d, "f_\eta"]
 K_\eta \arrow[l]  & \arrow[l, "p_\eta" above] K\arrow[d, "f"]\\
L_0 & \arrow[l, "q^\xi_0" above]  L_\xi & \arrow[l, "q^\eta_\xi" above]   L_\eta \arrow[l]  & \arrow[l, "q_\eta" above] L
\end{tikzcd}
\end{center}

Here $f_\xi$ is the restriction of $f$ to $K_\xi$, $L_\xi=f[K_\xi]$ and
the mapping $q_\xi$ is defined by \[ q_\xi(f(x))=f(p_\xi(x));\]
note that  property $(iii)$ guarantees that the definition is correct.
Then every $q_\xi$ is a retraction $L\to L_\xi$ and
again, for $\xi<\eta$ the bonding map $q^\eta_\xi$ is $q_\xi$ restricted to $L_\eta$.

For every $\xi<\lambda_0$ we have $w(K_\xi)\le |G_\xi|<\lambda$, so $L_\xi$ is $\kappa$-Corson
by the inductive assumption. Hence there is a family of functions $\Phi_\xi\sub C(L_\xi)$
which separates the points of $L_\xi$ and satisfies
$|\{f\in\Phi_\xi: f(x)\neq 0\}|<\kappa$ for every $x\in L_\xi$ (see Remark \ref{basic_observations}).
Consider
\[ \Phi=\{g\circ q_0:  g\in \Phi_0\}\cup
 \{g\circ q_{\xi+1} - g\circ q_\xi : \xi<\lambda_0, g\in \Phi_{\xi+1}\}\sub C(L);\]
 we shall check that $\Phi$ witnesses that $L$ is $\kappa$-Corson compact too.

Note first that $x=\lim_{\xi<\lambda} p_\xi(x)$, and therefore
\[ f(x)=\lim_{\xi<\lambda} f(p_\xi(x))= \lim_{\xi<\lambda} q_\xi(f(x))\]
 for every $x\in K$.

Take $x,y\in K$ such that $f(x)\neq f(y)$. If the points $f(x),f(y)$ are in $ L_0$, then some $g\in \Phi_0$
separates them. Otherwise, the first $\eta$ such that $q_\eta(f(x))\neq q_\eta(f(y))$ must be
a successor ordinal; if $\eta=\xi+1$, then
$f(x), f(y)$ are separated by $g\circ q_{\xi+1} - g\circ q_\xi$ for $g\in\Phi_{\xi+1}$ separating  $q_{\xi+1}(f(x)), q_{\xi+1}(f(y))$.

Given any $x \in K$, if
\[g\circ q_{\xi+1}(f(x)) - g\circ q_\xi(f(x))\neq 0\]
for $g\in \Phi_{\xi+1}$, then
\[ f(p_{\xi+1}(x))=q_{\xi+1}(f(x)) \neq   q_\xi(f(x))=f(p_\xi(x)), \]
so $p_{\xi+1}(x)\neq p_\xi(x)$.
We conclude that  $\supp(x)\cap (G_{\xi+1}\sm G_\xi)\neq\emptyset$, which may happen
$<\kappa$ many times.
Hence the cardinality of  $\vf\in\Phi$ with $\vf(f(x))\neq 0$ is $<\kappa$, and the proof is complete.
\end{proof}

Note that Theorem \ref{stable} says, in particular, that for regular $\kappa$,
every subalgebra of a $\kappa$-Corson algebra is again $\kappa$-Corson.
We do not know if the class of $\kappa$-Corson compacta is stable under continuous images also
for a singular $\kappa$. In particular, one can ask the following.

\begin{prob}
Is there a singular number  $\kappa$ and a subalgebra of a $\kappa$-Corson algebra which
is not $\kappa$-Corson?
\end{prob}

\section{Basically disconnected spaces and beyond}\label{bds}

We discuss   here  \cite[Question 5.2]{BKT}, a problem  which can be stated as follows.

\begin{prob}\label{p}
Can a $\sigma$-complete Boolean algebra $\fB$ of size $\kappa$ be
$\kappa$-Corson?
\end{prob}

Given an algebra $\fB$,  the statement `$\fB$ is $\sigma$-complete' can be weakened
in various directions, see Koszmider and Shelah \cite{KS13} and the references therein.
In particular, the following properties were considered in \cite{KS13}.

\begin{definition}
A Boolean algebra $\fB$ is said to have the  subsequential
completeness property (SCP) if for
 any sequence of pairwise disjoint $a_n\in\fB$ there is infinite
 $N\sub\omega$ such that the family $\{a_n:n\in N\}$ has a least upper bound in $\fB$.

A Boolean algebra $\fB$ is said to have the weak subsequential
separation property (WSSP) if for any sequence of pairwise disjoint $a_n\in\fB$
there is $b\in\fB$ such that
 both of the sets
$\{n \in\en : a_n\sub b\}$ and  $\{n \in\en : a_n\cap b=0\}$
are infinite.
\end{definition}

SCP was introduced by Haydon \cite{Ha81} in connection with the Grothendieck property of Banach spaces.
Clearly, we have the following implications (none of which is reversible, see  \cite{KS13})
\medskip

$\sigma$-complete $\longrightarrow$ SCP $\longrightarrow$   WSSP.
\medskip

Recall that a compact space $K$ is basically disconnected if the closure of
every  open $F_\sigma$ subset of $K$ is open.
In particular, the Stone space of a $\sigma$-complete Boolean algebra
is basically disconnected. In turn, such a property of a compact space $K$   can be weakened to saying
that $K$ is an $F$-space, i.e.\ every two disjoint open $F_\sigma$ subsets of $K$ have disjoint closures.
Recall that there are connected $F$-spaces,  e.g.\ $\beta([0,\infty))\setminus [0,\infty)$, see \cite[Theorem 14.27]{GJ60}.

Concerning Problem \ref{p},
we note first a couple of  easy answers.

\begin{rem}\label{ar} Let $\fB$ be  a Boolean algebra with $|\fB|=\kappa$.

\begin{enumerate}[(a)]
\item If   $\fB$ contains an independent family of size $\kappa$, then
$\ult(\fB)$ maps continuously onto $2^\kappa$, so $\ult(\fB)$ is not $\kappa$-Corson
compact by Theorem \ref{stable} and Lemma \ref{character}
(for regular $\kappa$; if $\kappa$ is singular, then use Lemma \ref{witek}).
\item In particular,  if  $\fB$ is complete, then it is not  $\kappa$-Corson, as
such  $\fB$ contains an independent family of full size by the Balcar-Franek theorem.
In the language of topology, no extremally disconnected compact space
of weight $\kappa$ is $\kappa$-Corson compact.
\item The answer to Problem \ref{p} is `no' for $\kappa=\cont$.
Indeed, taking pairwise disjoint $a_n\in\fB^+$ in a $\sigma$-complete algebra $\fB$ and
an independent family $\cI\sub P(\omega)$ of cardinality $\cont$ we get an independent family $\{a_I: I\in\cI\}$,
where
\[ a_I=\bigvee_{n\in I} a_n,\]
and we argue as above.
\item More generally, if $\fB$ has WSSP and  $\kappa=\cont$, then $\fB$ is not  $\cont$-Corson.
This is a consequence of \cite[Theorem 1.4]{KS13}, stating that such $\fB$ must contain an
independent family of size $\cont$.
\item
It is routine to check that
if $\fB$ has WSSP (in particular, if $\fB$ has SCP), then  $K=\ult(\fB)$
contains no nontrivial converging sequence.
This is a very particular version of \cite[Proposition 2.4]{KS13}
(related to the  Grothendieck property of the underlying space of continuous functions).
The lack of converging sequences and Lemma \ref{cs} might give a short proof of \ref{main} below
if we knew that the space $K$ in question would contain a point of character $\kappa$.
\end{enumerate}
\end{rem}

The following provides a negative answer to Problem \ref{p}
in case $\kappa$ in question is regular.

\begin{thm}\label{main}
If $\fB$ is a Boolean algebra of regular size $\kappa$ and $\fB$ has SCP, then
$\fB$ is not $\kappa$-Corson.
\end{thm}

\begin{proof} 
Suppose otherwise; we consider the Stone space $K=\ult(\fB)$ which is then
a zero-dimensional compact space of weight $\kappa$ and is $\kappa$-Corson compact
so, by Lemma \ref{lemma:ba},  we can assume that
\[K\sub\Sigma_\kappa(\{0,1\}^\kappa)=\{x\in \{0,1\}^\kappa: |\supp(x)|< \kappa\}.\]

For a  set $a\sub \kappa$ we write $C(a)=\{x\in K: x|a\equiv 1\}$;
if $\xi<\kappa$, then $C(\{\xi\})$ is simply denoted by $C(\xi)$.
Note that we can assume that $C(\xi)\neq\emptyset$ for every $\xi<\kappa$.

Using Lemma \ref{l2} we define an strictly increasing function $g:\kappa\to\kappa$ such that, for
every $\alpha<\kappa$,  the initial segment $\{\xi:\xi<g(\alpha)\}$ is good for $K$.
As an increasing union of good
sets is also good (Lemma \ref{la}), we can assume that $g$ is continuous, i.e.\
$g(\alpha)=\sup_{\beta<\alpha}g(\beta)$
whenever $\alpha$ is a limit ordinal. For simplicity, we shall write
\[ \pi_\alpha=\pi_{\{\xi:\xi<g(\alpha)\}}\mbox{ and }
 K_\alpha=\pi_{\alpha}[K]\sub K.\]

Consider now a family $\FF$ of sets $a\sub\kappa$ such that

\begin{enumerate}[(i)]
\item $a\sub g[\kappa]$,
\item  $C(a)\neq\emptyset$, and
\item $C(a\cup\{\xi\})=\emptyset$ for every $\xi>\max a$.
\end{enumerate}

\medskip

\noindent {\sc Claim 1.}
Every  $a\in \FF$ is finite.
\medskip

To verify Claim 1, suppose that $a\in\FF$ contains a strictly  increasing sequence of $\xi_n=g(\alpha_n)$;
set $\xi=g(\alpha)$, where $\alpha=\sup_n\alpha_n$.
Then there is $x\in K$ such that $x(\xi_n)=1$ for every $n$.
Now $\pi_{\alpha_n}(x)$ is a sequence of distinct points in $K$, and
$\lim_n \pi_{\alpha_n}(x)=\pi_\alpha(x)$, which contradicts the fact that $K$ contains no converging
sequences, see Remark \ref{ar}(e). 
\medskip

Next we consider a set $S\sub\kappa$ consisting of limit ordinals $\alpha<\kappa$
having the following property:
$\alpha=\sup_n\alpha_n$ for some $\alpha_n<\alpha$ and there are
$a_n\in\FF$ such that $a_n\sub g(\alpha_{n+1})\sm g(\alpha_n)$ for every $n$.
The following can be checked by a standard argument.
\medskip

\noindent {\sc Claim 2.}
The set $S$  is stationary.
\medskip

Note that, as $\fB$ has SCP,
 for every sequence of pairwise disjoint clopen sets $V_n\sub K$ there
is infinite $N\sub\omega$ such that
$ \overline{\bigcup_{n\in N} V_n}$
is a clopen subset of $K$.

For $\alpha\in S$ we pick a sequence of $a_n\in\FF$ witnessing that $\alpha$ belongs to $S$;
passing to a subsequence if necessary and applying the above remark, we
can assume that
\[ U_\alpha= \overline{\bigcup_{n\in N} C(a_n)}\]
is a clopen subset of $K$.

Applying  Lemma \ref{l4.5} we get the following.

\medskip

\noindent {\sc Claim 3.}
There is a function $f:S\to\kappa$ such that, for every $\alpha\in S$,  $f(\alpha)<\alpha$
and $U_\alpha\cap K_\alpha=C_\alpha\cap K_\alpha $, where $C_\alpha$ is a clopen set
determined by coordinates below $g(f(\alpha))$.
\medskip

By the pressing down lemma, $f(\alpha)=\eta<\kappa$ for $\alpha$ from some stationary set $S_0\sub S$.
It follows that for $\alpha\in S_0$, $C_\alpha$ is a clopen set determined by a finite number of
 coordinates below $g(\eta)<\kappa$.
 It follows that $C_\alpha=C$ for a fixed set $C$ and cofinally many $\alpha$'s.

We shall arrive at a contradiction considering $g(\eta)< \alpha<\beta$ in $S$ such that
\[ (*)\quad U_\alpha\cap K_\alpha=C\cap K_\alpha\mbox{ and } U_\beta\cap K_\beta=C\cap K_\beta.\]
Indeed, take two sequences of finite sets $a_n,b_n\in\cF$ witnessing that
$\alpha,\beta$ are in $S$, respectively.
Note that $J=\bigcup_n b_n \cap g(\alpha)$ is finite, so there is $k$ such that
$a_k\sub g(\alpha) \sm \max J$. Then $C(a_k)$ is disjoint from every $C(b_n)$, and hence
$C(a_k)\cap U_\beta=\emptyset$. Moreover, we have $C(a_k)\sub K_\alpha\sub K_\beta$, so
the second equation of (*) implies $C(a_k)\cap C=\emptyset$.

On the other hand, we have $C(a_k)\sub U_\alpha$, so $C(a_k)\sub C$ by the first equation of (*), and
this is in contradiction with $C(a_k)\neq\emptyset$.
This shows that (*) cannot hold, and the proof is complete.
\end{proof}

\begin{prob}
Does the assertion of Theorem \ref{main} remains true if we replace
SCP by WSSP?
\end{prob}

\section{Function spaces $C_p(K)$ on $\kappa$-Corson compacta $K$}\label{sec_CpK}

We briefly discuss in this section  generalizations of a
result due to Pol  \cite[\S 3, Thm.\ 1.1]{Po84} (see
also \cite{Gu77} and \cite{AP80}). Pol's theorem
characterizes  Corson compacta
$K$ in terms of topological properties of function spaces  $C_p(K)$. We extend this characterization for $\kappa$-Corson compacta, where $\kappa$ is regular, and use it to give an alternative proof of Theorem \ref{stable}.

Given  uncountable cardinal numbers $\lambda$ and $\kappa$, by $L_\kappa(\lambda)$
we denote the set $\lambda\cup\{\infty\}$ topologized as follows:
all points $\alpha\in\lambda$ are isolated in $L_\kappa(\lambda)$, and basic neighborhoods
of $\infty$ have the form $\{\infty\}\cup(\lambda\setminus A)$ where $A\sub\lambda$ and $|A|<\kappa$.

By $\elk$ we denote the class of all spaces which are continuous images of closed subsets of the countable product $L_\kappa(\lambda)^\omega$ for some cardinal $\lambda$.
Pol's result states   that  $K$ is Corson compact if and only if $C_p(K)\in \mathcal{L}_{\omega_1}$.
In \cite{BM04} this result was generalized for $\kappa^+$-Corson compact spaces.

The next proposition is a generalization of Proposition 4.1 from \cite{BM04} and can be justified in exactly the same way.

\begin{prp} \label{elk1} For every uncountable cardinal number $\kappa$, the class
$\elk$ is closed under taking
continuous images, closed subspaces, countable products and countable unions.
\end{prp}

Given an infinite cardinal number $\kappa$ and a topological space $X$, we say that $X$ is
$\kappa$-Lindel\"of if every open cover of $X$ contains a subcover of cardinality less than $\kappa$. In this terminology $\omega_1$-Lindel\"of spaces are usual Lindel\"of spaces.

Corollary 4.2 from Noble \cite{No71} implies that, for every uncountable cardinal $\lambda$, the space
$L_\kappa(\lambda)^\omega$ is $\kappa$-Lindel\"of, for every uncountable regular cardinal number $\kappa$, and it is $\kappa^+$-Lindel\"of, for every singular $\kappa$. Hence we obtain

\begin{prp} \label{elk2} Let $\kappa$ be an uncountable cardinal number and $X$ be an element of $\elk$. Then
\begin{enumerate}[(A)]
\item if $\kappa$ is regular, then $X$ is $\kappa$-Lindel\"of;
\item if $\kappa$ is singular, then $X$ is $\kappa^+$-Lindel\"of.
\end{enumerate}
\end{prp}

The next theorem can be proved using the same arguments as in the proof of one implication in Theorem 6.1 from \cite{BM04} (cf.\ \cite{Za})

\begin{thm}\label{thm_k_Corson_to_Lk}
Let $\kappa$ be an uncountable cardinal and $K$ be a $\kappa$-Corson compact space. Then the function space $C_p(K)$ belongs to $\elk$.
\end{thm}

\begin{cor}\label{cor_k_Corson_to_k_Lind}
Let $\kappa$ be an uncountable cardinal and $K$ be a $\kappa$-Corson compact space. Then
\begin{enumerate}[(A)]
\item if $\kappa$ is regular, then $C_p(K)$ is $\kappa$-Lindel\"of;
\item if $\kappa$ is singular, then $C_p(K)$ is $\kappa^+$-Lindel\"of.
\end{enumerate}
\end{cor}

The following fact belongs to the folklore

\begin{prp}\label{prp_ex_ord}
Let $\kappa$ be an infinite cardinal and $[0,\kappa]$ be the ordinal space equipped with the order topology. Then the function space $C_p([0,\kappa])$ contains a closed discrete subspace of cardinality $\kappa$.
\end{prp}

\begin{proof}
One can easily verify that the subspace
\begin{eqnarray*}
\{f\in C_p([0,\kappa]) \cap 2^{[0,\kappa]}: f(\kappa) = 0 \mbox{ and $f$ is nonincreasing}\} = \{\chi_{[0,\alpha]}: \alpha < \kappa\}\cup\{\chi_\emptyset\}
\end{eqnarray*}
has the required properties.
\end{proof}

Since the space $[0,\kappa]$ is $\kappa$-Corson for singular $\kappa$, and it is $\kappa^+$-Corson for regular $\kappa$, the above proposition shows that the results of Corollary \ref{cor_k_Corson_to_k_Lind} cannot be improved. In particular, this gives a negative answer to a part of Question 5.5 from \cite{BKT}.

The next theorem generalizes Theorem 6.1 from \cite{BM04}. One implication in this result is given by Theorem \ref{thm_k_Corson_to_Lk}. To prove the reverse implication one needs to modify the arguments from the proof of \cite[Theorem 6.1]{BM04} (which, in turn, follow the reasoning from \cite{Po84}). The details of the proof will be presented in \cite{Za}.

\begin{thm}\label{charact_k_Corson_CpK}
Let $\kappa$ be an uncountable regular  cardinal, and $K$ be a compact space. Then
$K$ is a $\kappa$-Corson compact space if and only if $C_p(K)\in\elk$.
\end{thm}

The next fact is well known

\begin{prp}\label{prp_surj_to_closed_emb}
If $L$ is a continuous image of a compact space $K$, then the function space $C_p(L)$ embeds as a closed subspace in $C_p(K)$.
\end{prp}

Using Theorem \ref{charact_k_Corson_CpK} and Propositions \ref{elk1} and \ref{prp_surj_to_closed_emb} we can immediately obtain an alternative proof of Theorem \ref{stable}.

\begin{prob}\label{prob_k_Corson_CpK}
Does Theorem \ref{charact_k_Corson_CpK} hold true for any uncountable cardinal number $\kappa$?
\end{prob}

\section{Measures on $\kappa$-Corson compacta}\label{measures}

In this section, we discuss the stability of the classes of $\kappa$-Corson compacta under the functor
$P$ which assigns to a compact space $K$ the compact space $P(K)$
of regular probability Borel measures equipped with the $weak^\ast$ topology.
Let us note that neither the class of $\omega$-Corson compacta nor $\ny$ is stable
in this sense, see Remark \ref{conv_no_NY}.

\begin{rem}\label{m:0}
Recall that the space $P(K)$ shares many topological properties with
$M_1(K)$, the space of all signed measures from the unit ball in $C(K)^\ast$.
In particular, $M_1(K)$ is a continuous image of
$T\times P(K)\times P(K)$, where $T=\{(t,s)\in\er^2: |t|+|s|\le 1\}$.
\end{rem}

If $P(K)$ is Corson compact, then so is $K$ since it is homeomorphic to
$\Delta_K=\{\delta_x:x\in K\}$. The reverse implication is not provable in ZFC.

 \begin{thm}\label{m:1}
The following are equivalent

\begin{enumerate}[(i)]
\item $P(K)$ is Corson compact for every Corson compact space $K$;
\item $M_1(K)$ is Corson compact for every Corson compact space $K$;
\item ${\rm cov} (2^{\omega_1})>\omega_1$.
\end{enumerate}
\end{thm}

Here ${\rm cov} (2^{\omega_1})>\omega_1$ means that the Cantor cube $2^{\omega_1}$
cannot be covered by
$\omega_1$ sets that are null with respect to the usual product measure.
Such a statement  follows from Martin's axiom MA$(\omega_1)$ and is simply denied by CH.

The equivalence (i)$\Leftrightarrow$(ii) is clear by \ref{m:0} while (i)$\Leftrightarrow$(iii)
is a combination of
results from \cite{AMN88} and \cite{KM95}.
Recall that under CH there are Corson compacta $K$ for which $P(K)$ is very far
from being Corson compact  ---  it can be derived from a result due to Talagrand \cite{Ta80}
that under CH there is a Corson compact space $K$ such that
$P(K)$ contains a copy of $\beta\omega$.

We show below that Theorem \ref{m:1} can
be naturally extended to $\kappa$-Corson compacta using
the following concept.

\begin{dfn}\label{m:2}
A regular cardinal $\kappa$ is a {\em caliber of Radon measures} if, for every
 compact space $K$, a measure $\mu\in P(K)$, and a family $\{A_{\alpha}:\alpha<\kappa\}$ in
 $Bor(K)$ of $\mu$-positive sets, there is $x\in K$ such that
 $\{\alpha<\kappa: x\in A_\alpha\}$ has cardinality $\kappa$.

\end{dfn}

Basic properties of calibers of measures are discussed in \cite{DP04}.
We recall here the most relevant facts.

\begin{thm}\label{m:3}
\begin{enumerate}[(a)]
\item A regular cardinal $\kappa$ is a caliber of Radon measures if and only if
it is a caliber of the usual product measure on $2^\kappa$.
\item The relation ${\rm cov} (2^{\omega_1})>\omega_1$ is equivalent to saying that
$\omega_1$ is a caliber of Radon measures.
\item $\cont^+$ is a caliber of Radon measures.
\end{enumerate}
\end{thm}

\begin{thm}\label{m:4} For a regular uncountable cardinal number $\kappa$ the following are equivalent

\begin{enumerate}[(i)]
\item $P(K)$ is $\kappa$-Corson compact for every $\kappa$-Corson compact space $K$;
\item $M_1(K)$ is $\kappa$-Corson compact for every $\kappa$-Corson compact space $K$;
\item For every $K\sub\sikg$ and $\mu\in P(K)$ there is $S\sub\Gamma$ such that $|S|<\kappa$ and
\[\mu\left(\{x\in K; \supp(x)\sub S\}\right)=1;\]
\item $\kappa$ is a caliber of Radon measures.
\end{enumerate}

\end{thm}

\begin{proof}
The equivalence (i)$\Leftrightarrow$(ii) follows from $P(K)\sub M_1(K)$, Theorem \ref{stable}
and Remark \ref{m:0}.

For (i)$\Rightarrow$(iii) consider $K\sub \sikl$ and $\mu\in P(K)$. It is not difficult to check that
\[ \mu\in \overline{{\rm conv} \{\delta_x:x\in K\} }.\]
By Lemma \ref{st} applied to $P(K)$, there is a set $D\sub K$ such that $|D|<\kappa$ and
\[ \mu\in \overline{{\rm conv} \{\delta_x:x\in D\} }.\]
We shall check that  the following set satisfies (iii):
\[ S=\bigcup_{x\in D} \supp(x).\]

 Write $C(\gamma)=\{x\in K: x(\gamma)=0\}$; if $\gamma\notin S$
then $\delta_x(C(\gamma))=1$ for every $x\in D$ and hence
$\nu(C(\gamma))=1$ for every $\nu$ from the  set $ {\rm conv} \{\delta_x:x\in D\}$.
Since $\mu$ is in its closure, we have $\mu(C(\gamma))=1$ as well.
Finally, the set
\[ \{x\in K: \supp(x)\sub S\}=\bigcap_{\gamma\in \Gamma\sm S} C(\gamma)\]
is an intersection of closed sets of full measure $\mu$ so it has measure 1, by regularity of  $\mu$.

To check (iii)$\Rightarrow$(i) take  $K\sub \sikg$ and consider
the family $\cF \sub C(K)$ of functions $f_{(I, \vf)}$ where
\[ f_{(I, \vf)} (x)=\prod_{\gamma\in I} x(\gamma)^{\vf(\gamma)},
\quad I\sub [\Gamma]^{<\omega}, \vf:I\to\en.\]
 By the Stone-Weierstrass theorem,
the linear span of $\cF$ is norm dense in $C(K)$, and therefore the family $\cF$
distinguishes elements of $P(K)$. Given $\mu\in P(K)$, take a set $S$ as in $(iii)$;
then $\mu(f_{(I, \vf)})>0$ implies $I\sub S$ so, treating $\cF$
as a family of  continuous functions on $P(K)$ and using Remark \ref{basic_observations}(d),
we conclude that $P(K)$ is $\kappa$-Corson.

To argue for (iii)$\Rightarrow$(iv), 
suppose that $\kappa$ is not a caliber of a measure $\mu\in P(L)$
for some compact space $L$. This is witnessed, in view of  regularity of $\mu$, by
a family $\cA=\{A_\xi: \xi<\kappa\}$ of closed subsets of $L$ such that $\mu(A_\xi)>0$ for every $\xi$ and
\[ \Big| \{\xi<\kappa: y\in A_\xi\}\Big|<\kappa,\]
for every $y\in L$. Consider the Boolean algebra $\fB$ of subsets of $L$   generated by $\cA$.
Then $\fB$ is a $\kappa$-Corson algebra
by Lemma \ref{lemma:ba}.
Moreover,
\[ \theta:\ult(\fB)\to 2^\cA,   \theta (p)=\chi_{\widehat{A}}(p) \mbox{ for } p\in\ult(\fB),\]
is an embedding onto $K\sub \Sigma_\kappa(2^\cA)$. The measure $\mu$ defines
a measure $\widehat{\mu}\in P(\ult(\fB))$ via the Stone isomorphism, and if we take
the image measure $\nu=\theta[\widehat{\mu}]\in P(K)$, then
\[ \nu(\{x\in K: x(A)=1\}=\mu(A)>0,\]
for every $A\in\cA$, so we have proved that (iii) does not hold.

Finally, we verify (iv)$\Rightarrow$(iii). If $K\sub\sikg$ and $\mu\in P(K)$, then the family
of sets of the form $C_\gamma=\{x\in K: x(\gamma)>0\}$, $\gamma\in\Gamma$, does
not contain $\kappa$ many sets with nonempty intersection so, by $(iv)$,
the set $S=\{\gamma\in\Gamma: \mu(C_\gamma)>0\}$ must be of size $<\kappa$
and hence  $\mu$ is supported by $\{x\in K: \supp(x)\sub S\}$.
\end{proof}

\begin{cor}\label{m:5}
For every $\cont^+$-Corson compact space $K$ the spaces
$P(K)$ and $M_1(K)$ are $\cont^+$-Corson compact.
\end{cor}

\section{On Banach spaces $C(K)$}\label{banach}

Let us say that a class of compact spaces $\cC$ is B-stable if for any $L\in\cC$ and
any compact space $K$, the fact that the Banach spaces $C(K)$ and $C(L)$ are isomorphic implies
$K\in\cC$. Clearly, the class of metrizable compacta is B-stable, since
$K$ is metrizable if and only if $C(K)$ is separable.
 The class of Eberlein compacta is B-stable, as $L$ is Eberlein compact if and only if
 $C(L)$ is weakly compactly generated, which is
an isomorphic property, see e.g.\ \cite{Ne}.

The class of $\omega$-Corson compact spaces $K$ is not B-stable:
$C(I)$ and $C(I^\omega)$ are isomorphic by Miljutin's theorem (cf.\ Semadeni \cite{Se}),
while $I^\omega$ is not $\omega$-Corson.
We do not know if the class $\ny$ is B-stable.
The following open question has been around for several years,
see e.g. \cite[3.9]{AMN88} and \cite[6.45]{Ne}.

\begin{prob}\label{bs:1}
Is the class of Corson compacta B-stable?
\end{prob}

The problem has a relatively easy  answer `yes' if  $\omega_1$ is
a caliber of Radon measures. However, it is still unclear if \ref{bs:1} can be resolved within the usual axioms
of set theory; there are partial positive solutions in \cite{Pl15i}.
It seems natural to  formulate a more general problem.

\begin{prob}\label{bs:2}
Given an uncountable  cardinal number $\kappa$, is the class of $\kappa$-Corson compacta  B-stable?
\end{prob}

Again, the positive answer to \ref{bs:2} follows for
every  regular number $\kappa$ for which the class of $\kappa$-Corson compacta
is stable under $P(\cdot)$.
In fact, we have the following more general result.

\begin{thm}\label{bs:3}
Suppose that a regular cardinal number  $\kappa$ is a caliber of Radon measures.
If $K$ is $\kappa$-Corson compact and $T: C(L)\to C(K)$ is an isomorphic embedding, then
$L$ is $\kappa$-Corson compact.
\end{thm}

\begin{proof}
We can assume that
$c\cdot \|g\|\le \|Tg\|\le \|g\| $ for some $c>0$ and every $g\in C(L)$.
We consider the dual operator
\[T^\ast: M(K)\to M(L),\quad  T^\ast\mu(g)=\mu(Tg)
\mbox{ for } \mu\in M(K) \mbox{ and }g\in C(L).\]
Then $T^\ast$ is surjective and $weak^\ast$ -- $weak^\ast$
continuous; moreover,  $T^\ast[M_1(K)]\supseteq c\cdot M_1[L]$.
By Theorem \ref{m:4}, $M_1[K]$ is $\kappa$-Corson compact, so
$T^\ast[M_1(K)]$ is also $\kappa$-Corson compact.
The space $L$ embeds into $T^\ast[M_1(K)]$ via the mapping $L\ni y\mapsto c\cdot\delta_y$,
and we are done.
\end{proof}

 Corollary \ref{m:5} and Theorem \ref{bs:3} yield the following.

\begin{cor}
The class of $\cont^+$-Corson compacta is B-stable.
\end{cor}


\begin{thebibliography}{AMN}
\bibitem[Al]{Al79}
K.\  Alster, {\em  Some remarks on Eberlein compacts},
Fund.\ Math.\ 104 (1979),  43--46.

\bibitem[AP]{AP80}
K.\ Alster and R.\ Pol, \emph{On function spaces of compact
subspaces of $\Sigma$-products of the real line}, Fund. Math. {107}
(1980), 135--143.

\bibitem[AMN]{AMN88}
S.\  Argyros, S.\  Mercourakis and S.\ Negrepontis,
{\em Functional-analytic properties of Corson-compact spaces},  Studia Math.\ 89 (1988),  197--229.	

\bibitem[BJ]{BJ} T.\ Bartoszyński and H.\ Judah,
\emph{Set Theory. On the Structure of the Real Line},  A K Peters, Ltd., Wellesley, MA, (1995).

\bibitem[BM]{BM04}
M.\ Bell and W.\  Marciszewski, {\em  Function spaces on $\tau$-Corson compacta and tightness of
polyadic spaces},  Czechoslovak Math. J. 54(129) (2004), 899--914.
	
\bibitem[BKT]{BKT} R.\ Bonnet, W.\ Kubi\'s and S.\ Todor\v{c}evi\'c,
\emph{Ultrafilter selection and Corson compacta},
Rev.\ R.\ Acad.\ Cienc.\ Exactas Fís.\ Nat.\ Ser.\ A Mat.\ RACSAM 116 (2022),  Paper No.\ 178, 26 pp.	

\bibitem[vD]{vD} E.\ van Douwen,  \emph{The integers and topology, Handbook of Set-Theoretic Topology} (K. Kunen
	and J. Vaughan, eds.), North-Holland, 1984, pp.\ 111--167.
	
\bibitem[DP]{DP} A.\ Dow and E.\ Pearl, \emph{Homogeneity in powers of zero-dimensional first-countable spaces}, Proc. Amer. Math. Soc. \textbf{125} (1997), 2503--2510.

\bibitem[DPl]{DP04}
 M.\ D\v{z}amonja and G.\  Plebanek,
{\em  Precalibre pairs of measure algebras},  Topology Appl. 144 (2004),  67--94.
	
\bibitem[En1]{En1} 	R.\ Engelking, \emph{General Topology}, Heldermann Verlag, Berlin (1989).
	
	
\bibitem[En2]{En}	R.\ Engelking, {\em Theory of Dimensions Finite and Infinite}, Sigma Series in Pure Mathematics, 10. Heldermann Verlag, Lemgo, 1995.
	
\bibitem[EP]{EP} R.\ Engelking and R.\ Pol, {\em Compactifications of countable-dimensional and strongly countable-dimensional spaces},
	Proc. Amer. Math. Soc. \textbf{104} (1988), 985--987.

\bibitem[GJ]{GJ60} L. Gillman and M. Jerison, \textit{Rings of continuous functions}, The University Series in
Higher Mathematics, D. Van Nostrand Co., Inc., Princeton, N.J.--Toronto--London--New York,
1960.
	
\bibitem[Gr]{Gr} G.\ Gruenhage,
	\emph{Covering properties on $X^2\setminus \Delta$, W-sets, and compact subsets of $\Sigma$-products},
	Topology Appl. \textbf{17} (1984), no. 3, 287--304.

\bibitem[Gu]{Gu77} S.P.\ Gul'ko, {\em
Properties of sets that lie in $\Sigma $-products},
Dokl.\ Akad.\ Nauk SSSR 237 (1977),  505--508.

\bibitem[Ha]{Ha81}
R.\ Haydon, {\em  A nonreflexive Grothendieck space that does not contain $\ell_\infty$},
Israel J.\ Math.\  40 (1981), 65–-73.

\bibitem[Ka1]{Ka99}
O.\ Kalenda,  {\em
Embedding of the ordinal segment $[0,\omega_1]$ into continuous images of Valdivia compacta},
Comment.\ Math.\ Univ.\ Carolin.\ 40 (1999),  777--783.
	
\bibitem[Ka2]{Ka} O.\ Kalenda, {\em  Valdivia compact spaces in topology and Banach space theory}, Extracta Math. \textbf{15} (2000),
	1–-85.

\bibitem[KS]{KS13}
 P.\  Koszmider and S.\  Shelah, {\em Independent families in Boolean algebras with some separation properties},
  Algebra Universalis 69 (2013),  305--312.


\bibitem[KvM]{KM95}
 K.\ Kunen and  J.\  van Mill, {\em Measures on Corson compact spaces}, Fund.\ Math.\ 147
(1995), 61--72.
\bibitem[MP]{MP17}
M.\ Magidor and G.\ Plebanek, {\em  On properties of compacta that do not reflect in small continuous images},
Topology Appl.\ 220 (2017), 131--139.	
\bibitem[Ma]{Ma} W.\ Marciszewski, \emph{On two problems concerning Eberlein compacta}, Rev. R. Acad. Cienc. Exactas Fís. Nat. Ser. A Mat. RACSAM \textbf{115} (2021), no. 3, Paper No. 126, 14 pp.

\bibitem[MR]{MR77}
E.\ Michael and  M.E.\  Rudin, {\em A note on Eberlein compacts},
Pacific J.\ Math.\ 72 (1977),  487--495.
	
\bibitem[NY]{NY} L.B.\ Nakhmanson and N.N.\  Yakovlev,
	\emph{Bicompacta lying in $\sigma$-products},
	Comment. Math. Univ. Carolin. \textbf{22} (1981), no. 4, 705--719.
	
\bibitem[Ne]{Ne} S.\ Negrepontis, \emph{Banach spaces and topology,
		Handbook of Set-Theoretic Topology} (K. Kunen
	and J. Vaughan, eds.), North-Holland, 1984, pp.\  1045--1142.

\bibitem[No]{No71}
N.\ Noble, \emph{Products with closed projections. II}, Trans. Amer. Math.
Soc. {160} (1971), 169--183.

\bibitem[Ny1]{Ny1} P.J.\ Nyikos, \emph{On the product of metacompact spaces I. Connections with hereditary compactness},  Amer.\ J.\ Math.\ 100 (1978), no. 4, 829--835.
	
\bibitem[Ny2]{Ny2} P.J.\ Nyikos, \emph{Various topologies on trees} Proceedings of the Tennessee Topology Conference, World Scientific Publishing Co., 1997, pp. 167--198.


\bibitem[Pl1]{Pl15i} G.\ Plebanek,
{\em On isomorphisms of Banach spaces of continuous functions},
 Israel J.\ Math.\ 209 (2015),  1--13.

\bibitem[Pl2]{Pl22} G.\ Plebanek,
{\em Musing on Kunen's compact $L$-space}, Topology Appl.
 323 (2023), Paper No. 108294, 10 pp.

\bibitem[Po]{Po84}
R.\  Pol,  {\em On weak and pointwise topology in function spaces},
University of Warsaw preprint No. 4184, Warsaw, 1984.
	
\bibitem[Se]{Se}
	Z.\ Semadeni, \textit{Banach Spaces of Continuous Functions}, PWN, Warsaw, 1971.

\bibitem[Ta]{Ta80}
M.\  Talagrand, {\em Sur les mesures vectorielles d\'efinies par une application Pettis-int\'egrable},
 Bull.\ Soc.\ Math.\ France 108 (1980), no. 4, 475--483.

\bibitem[Tk]{Tk83} M.G.\ Tka\v{c}enko, {\em Sum theorems for the tightness and $\pi$character},
Comm.\ Math.\ Univ.\ Car.\ 24 (1983), 51--62.

	
\bibitem[WZ]{Order} S.W.\ Williams and H.\ Zhou, \emph{The order-like structure of compact monotonically normal spaces}, Comment.\ Math.\ Univ.\ Carolin. 39 (1998), 207–-217

	
\bibitem[Ya]{Ya} N.N.\ Yakovlev,
	\emph{On bicompacta in $\Sigma$-products and related spaces},
Comment.\ Math.\ Univ.\ Carolin.\ {21} (1980),  263--283.

\bibitem[Za]{Za} K.\ Zakrzewski, \emph{Function spaces on Corson-like compacta}, in preparation.
	
\end{thebibliography}
\end{document}